\documentclass[11pt]{article}
\usepackage{amsmath,amssymb,mathrsfs}
\usepackage{amsfonts}
\usepackage{color}
\newtheorem{theo}{Theorem}[section]
\newtheorem{lem}[theo]{Lemma}
\newtheorem{cor}[theo]{Corollary}

\newcommand{\mysection}[1]{\section{#1} \setcounter{equation}{0}}
\newcommand{\proof}{{\sc Proof.} \quad}
\newcommand{\proofc}{{\sc Proof} \ }
\newcommand{\be}{\begin{equation} \label}
\newcommand{\ee}{\end{equation}}
\newcommand{\bea}{\begin{eqnarray}\label}
\newcommand{\eea}{\end{eqnarray}}
\newcommand{\bas}{\begin{eqnarray*}}
\newcommand{\eas}{\end{eqnarray*}}
\newcommand{\bit}{\begin{itemize}}
\newcommand{\eit}{\end{itemize}}
\newcommand{\qed}{\hfill$\Box$ \vskip.2cm}
\newcommand{\nn}{\nonumber}
\newcommand{\R}{\mathbb{R}}
\newcommand{\N}{\mathbb{N}}

\newcommand{\eps}{\varepsilon}

\newcommand{\abs}{\\[5pt]}

\newcommand{\ou}{\overline{u}}

\newcommand{\uz}{\underline{z}}

\newcommand{\tu}{\widetilde{u}}
\newcommand{\tl}{\widetilde{\cal L}}
\newcommand{\os}{\overline{s}}
\begin{document}
\title{A Gagliardo-Nirenberg-type inequality and its applications \\
to decay estimates for solutions of a degenerate parabolic equation}
\author{
Marek Fila\footnote{fila@fmph.uniba.sk, corresponding author}\\
{\small Department of Applied Mathematics and Statistics, Comenius University,}\\
{\small 84248 Bratislava, Slovakia}
\and
Michael Winkler\footnote{michael.winkler@math.uni-paderborn.de}\\
{\small Institut f\"ur Mathematik, Universit\"at Paderborn,}\\
{\small 33098 Paderborn, Germany} }
\date{}
\maketitle
\begin{abstract}
\noindent 
We establish a Gagliardo-Nirenberg-type inequality in $\R^n$ for functions
which decay fast as $|x|\to\infty$. We use this inequality to derive upper
bounds for the decay rates of solutions of a degenerate parabolic equation.
Moreover, we show that these upper bounds, hence also the
Gagliardo-Nirenberg-type inequality, are sharp in an appropriate sense.
\abs
\end{abstract}
\newpage
\mysection{Introduction}
{\bf A Gagliardo-Nirenberg-type inequality.}\quad
The Gagliardo-Nirenberg inequalities (\cite{Gag1}, \cite{Gag2}, \cite{Nir}) play an important role in studying partial differential equations (cf. \cite{BB}, for instance). Consider the special case when $1\leq r<q<\infty$, 
and $0\leq\theta\leq 1$ are such that
\[
\frac{1}{q}=\frac{\theta}{r}+(1-\theta)\left(\frac{1}{2}-\frac{1}{n}\right).
\]
Then there is a constant $c>0$ which depends only on $n, q$ and $r$, such that 
any $\varphi\in L^r(\R^n)$ with $\nabla\varphi\in L^2(\R^n)$ satisfies
\be{GN}
	\|\varphi\|_{L^q(\R^n)} \le c \|\varphi\|_{L^r(\R^n)}^\theta 
	\|\nabla\varphi\|_{L^2(\R^n)}^{1-\theta}.
\ee
Our aim is to establish a new optimal inequality of a similar type by replacing the term $\|\nabla\varphi\|_{L^2(\R^n)}^{1-\theta}$ with $\|\nabla\varphi\|_{L^2(\R^n)}F\left(\|\nabla\varphi\|_{L^2(\R^n)}\right)$, 
where $F$ is some positive function with the properties that
\bas
	F(s)\to\infty
	\mbox{ as } s\to 0
	\qquad \mbox{and} \qquad
 	s^\theta F(s)\to 0
	\mbox{ as } s\to 0 \mbox{ for any $\theta>0$.}
\eas
The term 
$c \|\varphi\|_{L^r(\R^n)}^\theta$ will then be replaced by a constant which depends only on $n, q$ and $\int_{\R^n}{\cal L}(\varphi)$, where ${\cal L}$ is a suitable function related to $F$. The integrability of ${\cal L}(\varphi)$ will require fast (exponential-like) decay of $\varphi$.\abs
The  Gagliardo-Nirenberg inequalities (GNI) from \cite{Gag1}, \cite{Gag2}, \cite{Nir} have been improved and 
extended in many different directions. We shall mention some examples below without trying to give an exhaustive list.
Sharp constants in GNI in $\R^n$ were studied in \cite{Ag1, Ag2, Ag3, CENV,
dPD, DELL, DFLP, DT1, KSW, KW, LYZ, Zhai} and in GNI on Riemannian manifolds in
\cite{ACM, CM1, CM2}. The sharp constant in an anisotropic GNI with
fractional derivatives was found in \cite{Esf}. A pointwise GNI can be found
in \cite{MSh}, a weighted GNI in \cite{DV} and a GNI on manifolds in
\cite{Badr}. GNI in Orlicz spaces were
established in \cite{KK, KPP1, KPP2}, in Besov spaces of negative order in
\cite{MM}, in weak Lebesgue spaces in \cite{MCRR} and in spaces of functions
with bounded mean oscillation in \cite{MCRR, RS}. An affine GNI was derived
in \cite{LYZ, Zhai} and a nonlinear GNI in \cite{RS}.  
Connections between logarithmic Sobolev inequalities and generalizations of GNI were investigated in \cite{BG}.\abs
We are not aware, however, of any example of a GNI in the literature which to an essentially optimal extent
makes use of a presupposed superalgebraically fast decay of the involved function. 
Addressing this problem, as the main result of this paper we shall obtain the following. 
\begin{theo}\label{theo11}
  Assume that $s_0>0$ and ${\cal L}\in C^0([0,\infty)) \cap C^1((0,s_0))$ is positive on $(0,\infty)$, bounded,
  nondecreasing and such that the following holds:
  \bit
  \item[{\rm (H)}] There exist $a, \lambda_0>0$ such that
  \be{4.1}
        {\cal L}(s) \le (1+a\lambda) {\cal L}(s^{1+\lambda})
        \qquad \mbox{for all $s\in (0,s_0)$ and } \lambda\in (0,\lambda_0).     
  \ee 
  \eit
  Then for any $K>0$ and $q>0$ such that $q<\frac{2n}{(n-2)_+}$
  there exists $C=C(n,q,K)>0$ such that if
  $0\not\equiv \varphi\in W^{1,2}(\R^n)$ is a nonnegative function satisfying
  \be{11.1}
	\int_{\R^n} {\cal L}(\varphi) \le K,
  \ee
  then
  \be{11.3}
	\|\varphi\|_{L^q(\R^n)} \le C \|\nabla\varphi\|_{L^2(\R^n)} 
	{\cal L}^{-(\frac{1}{q}-\frac{n-2}{2n})} \Big(\|\nabla\varphi\|_{L^2(\R^n)}^2 \Big).
  \ee
\end{theo}
Typical examples of functions fulfilling (H) are given by
\be{el}
{\cal L}(s)=\ln^{-\kappa}\frac{M}{s} \qquad\mbox{ or }\qquad
{\cal L}(s)=\ln^{-\kappa}\ln\frac{M}{s} \, ,\qquad s>0, \kappa>0, M>e,
\ee
see Lemmata~\ref{lem15} and \ref{lem21}.
With such a choice of ${\cal L}$, (\ref{11.1}) will be satisfied if there are positive constants $c_0, \alpha, \beta$ and $\gamma$ such that
\[
\varphi(x) \le c_0 \, e^{-\alpha |x|^\beta}\qquad\mbox{ or }\qquad 
\varphi(x) \le c_0 \, \exp \Big\{ -\alpha \exp \big(\beta|x|^\gamma\big)\Big\}
\qquad \mbox{for all } x\in\R^n,
\]
respectively.
Notice that if $n>2$, $r:=\frac{n}{n-2}$ and ${\cal L}(s)=s^r$, then (\ref{11.3}) with
\[
C=c \|\varphi\|_{L^r(\R^n)}^\theta,\qquad \theta=\frac{2n}{q(n-2)} -1\, ,
\]
corresponds to (\ref{GN}).\abs
We shall next show that the exponent $\frac{1}{q}-\frac{n-2}{2n}$ in (\ref{11.3}) is sharp. 
We accomplish that in the context of an analysis of temporal decay rates in a degenerate parabolic equation,
for which Theorem~\ref{theo11} will yield certain upper bounds that thereafter, essentially 
by means of arguments based on parabolic comparison principles, will be seen to be optimal
in an appropriate sense.

\bigskip
{\bf Applications to decay estimates for a degenerate parabolic equation.}\quad
For $p\ge 1$,
consider the Cauchy problem
\be{0}
	\left\{ \begin{array}{l}
	u_t=u^p \Delta u, \qquad x\in\R^n, \ t>0, \\[1mm]
	u(x,0)=u_0(x), \qquad x\in\R^n,
 	\end{array} \right.
\ee
where $u_0 \in C^0(\R^n)\cap L^\infty(\R^n)$.
Our purpose is to study the large time behavior of global classical solutions under the assumption that
\be{growth}
	u_0(x)\to 0	\qquad \mbox{as } |x|\to\infty,
\ee
and our particular focus is on describing in a quantitative manner how various types of decay of $u_0$
affect the asymptotic behavior of $\|u(\cdot,t)\|_{L^\infty(\R^n)}$ as
$t\to\infty$.\abs
Before addressing this, as a caveat we need to note that even in the framework of smooth positive solutions, 
uniqueness does not hold for (\ref{0}). 
After all, however, (\ref{0}) always possesses a {\em minimal} global classical solution
$u$ for any positive continuous and bounded initial data (\cite{fast_growth1}).
This solution is minimal in the sense that whenever $T\in (0,\infty]$ and
$\tu \in C^0(\R^n\times [0,T)) \cap C^{2,1}(\R^n\times (0,T))$ are such that $\tu$ is positive and
solves (\ref{0}) classically in $\R^n\times (0,T)$ then we have $u\le \tu$ in $\R^n\times (0,T)$.\abs
Now for any initial data decaying sufficiently fast in space, this minimal solution is known to approach
zero at a temporal rate which at its leading order is determined by the algebraic function $t^{-\frac{1}{p}}$,
but which in fact must involve a subalgebraic correction. 
More precisely, the following was shown in \cite{fast_growth1}.
\begin{theo}\label{theo200}
  If $p\ge 1$ and $u_0\in \bigcap_{q_0>0} L^{q_0}(\R^n)$, then for any $\delta>0$ one can
  find
  $C(\delta)>0$ such that for the minimal solution $u$ of {\rm (\ref{0})} we have
  \be{200.3}
	\|u(\cdot,t)\|_{L^\infty(\R^n)} \le C(\delta) t^{-\frac{1}{p}+\delta}
	\qquad \mbox{for all } t>0.
  \ee
  Any global positive classical $u$ of {\rm (\ref{0})}
  has the property that for every $R>0$,
  \bas
        \inf_{|x|<R} \Big\{ t^\frac{1}{p} u(x,t) \Big\} \to +\infty
        \qquad \mbox{as } t\to\infty.
  \eas
\end{theo}

This theorem suggests that logarithmic terms may occur in the sharp decay
rate of $\|u(\cdot,t)\|_{L^\infty(\R^n)}$ if the decay of $u_0$ is fast
enough. We show that Theorem~\ref{theo11} implies an upper bound which
supports this conjecture.
\begin{theo}\label{theo14}
  Suppose that $p\ge 1$, that $s_0>0$, and that 
  ${\cal L}\in C^0([0,\infty)) \cap C^2((0,s_0))$ is positive and nondecreasing on $(0,\infty)$ with ${\cal L}(0)=0$ and 
  such that {\rm (H)} is valid, and such that furthermore
  \be{14.1}
	s{\cal L}''(s) \ge - \frac{3p+q_0-2}{p+q_0} {\cal L}'(s)
	\qquad \mbox{for all } s\in (0,s_0)
  \ee
  with a certain $q_0>0$.
  Moreover, assume 
  that $u_0\in C^0(\R^n)$ is positive, radially symmetric and nonincreasing with
  respect to $|x|$ and such that
  \be{14.01}
	u_0<\min \Big\{ s_0^\frac{2}{p} \, , \, s_0^\frac{2}{p+q_0} \Big\}
	\qquad \mbox{in } \R^n
  \ee
  as well as
  \be{14.2}
	\int_{\R^n} {\cal L}(u_0) < \infty.
  \ee
  Then there exist $t_0>0$ and $C>0$ such that the minimal solution $u$ of
  {\rm (\ref{0})} satisfies
  \be{14.3}
	\|u(\cdot,t)\|_{L^\infty(\R^n)} \le C t^{-\frac{1}{p}} {\cal L}^{-\frac{2}{np}} \Big(\frac{1}{t}\Big)
	\qquad \mbox{for all } t\ge t_0.
  \ee
\end{theo}
%
%
%
%
%
%
%
As observed in Lemmata~\ref{lem15} and \ref{lem21} below, 
the condition (\ref{14.1}) is indeed satisfied by the functions from (\ref{el}).
Firstly, concentrating on the first choice therein, as a consequence of Theorem~\ref{theo14} we obtain the following.
\begin{cor}\label{cor17}
  Let $p\ge 1$, and suppose that $u_0\in C^0(\R^n)$ be positive and such that
  \be{17.1}
	u_0(x) \le c_0 \, e^{-\alpha |x|^\beta}
	\qquad \mbox{for all } x\in\R^n
  \ee
  with positive constants $c_0,\alpha$ and $\beta$.
  Then for any $\delta>0$ one can find $t_0>0$ and $C>0$ such that the minimal solution of
  {\rm (\ref{0})} satisfies
  \be{17.2}
	\|u(\cdot,t)\|_{L^\infty(\R^n)}
	\le C t^{-\frac{1}{p}} \ln^{\frac{2}{p\beta}+\delta} t
	\qquad \mbox{for all } t\ge t_0.
  \ee
\end{cor}
This refines Theorem~\ref{theo200} and we shall see below that the upper
bound (\ref{17.2}) is optimal in an appropriate sense.
The second option offered by (\ref{el}) indicates that
for initial data with faster decay, 
also iterated logarithms may occur in the upper bounds:
%
%
%
%
%
%
%
%
%
\begin{cor}\label{cor23}
  Let $p\ge 1$, and assume that $u_0\in C^0(\R^n)$ is positive and such that there exist positive constants
  $c_0,\alpha,\beta$ and $\gamma$ fulfilling
  \be{23.1}
	u_0(x) \le c_0 \, \exp \Big\{ -\alpha \exp \big(\beta|x|^\gamma\big)\Big\}
	\qquad \mbox{for all } x\in\R^n.
  \ee
  Then for all $\delta>0$ one can find $t_0>e$ and $C>0$ such that the minimal solution of
  {\rm (\ref{0})} satisfies
  \be{23.2}
	\|u(\cdot,t)\|_{L^\infty(\R^n)}
	\le C t^{-\frac{1}{p}} \ln^{\frac{2}{p\gamma}+\delta} \ln t
	\qquad \mbox{for all } t\ge t_0.
  \ee
\end{cor}
From a corresponding lower bound below one can see that (\ref{23.2}) is also
sharp. This lower bound will follow from our next result. 
\begin{theo}\label{theo31}
  Let $p\ge 1$, and let $\Lambda\in C^0([0,\infty))$ be strictly increasing and such that
  \be{31.1}
	\frac{\Lambda(s)}{\ln s} \to +\infty
	\qquad \mbox{as } s\to\infty.
  \ee
  Moreover, assume that $u$ is a positive classical solution of {\rm (\ref{0})} in $\R^n\times (0,\infty)$, with initial data
  $u_0\in C^0(\R^n)$ satisfying
  \be{31.2}
	u_0(x)\ge e^{-\Lambda(x)}
	\qquad \mbox{for all } x\in\R^n.
  \ee
  Then there exist $t_0>0$ and $C>0$ such that
  \be{31.3}
	\|u(\cdot,t)\|_{L^\infty(\R^n)} 
	\ge C t^{-\frac{1}{p}} \Big\{ \Lambda^{-1} \big( C\ln t\big) \Big\}^\frac{2}{p}
	\qquad \mbox{for all } t\ge t_0,
  \ee
  where $\Lambda^{-1}$ denotes the inverse of $\Lambda$.
\end{theo}
For particular choices of $\Lambda$ we have the following two consequences:
\begin{cor}\label{cor32}
  Let $p\ge 1$, and suppose that $u_0\in C^0(\R^n)$ is such that there exist $c_0, \alpha, \beta>0$ fulfilling
  \be{32.1}
	u_0(x) \ge c_0 \, e^{-\alpha |x|^\beta}
	\qquad \mbox{for all } x\in \R^n.
  \ee
  Then one can find $t_0>1$ and $C>0$ such that any positive classical solution of
  {\rm (\ref{0})} satisfies
  \be{32.2}
	\|u(\cdot,t)\|_{L^\infty(\R^n)}
	\ge C t^{-\frac{1}{p}} \ln^\frac{2}{p\beta} t
	\qquad \mbox{for all } t\ge t_0.
  \ee
\end{cor}
\begin{cor}\label{cor33}
  Let $p\ge 1$ and $u_0\in C^0(\R^n)$ be such that
  \be{33.1}
	u_0(x) \ge c_0 \, \exp \Big\{-\alpha \exp (\beta|x|^\gamma)\Big\}
	\qquad \mbox{for all } x\in\R^n
  \ee
  with positive constants $c_0,\alpha,\beta$ and $\gamma$.
  Then there exist $t_0>e$ and $C>0$ with the property that any positive classical solution of
  {\rm (\ref{0})} satisfies
  \be{33.2}
	\|u(\cdot,t)\|_{L^\infty(\R^n)}
	\ge C t^{-\frac{1}{p}} \ln^\frac{2}{p\gamma} \ln t
	\qquad \mbox{for all } t\ge t_0.
  \ee
\end{cor}
These last two corollaries imply that the upper bounds (\ref{17.2}) and
(\ref{23.2}) cannot hold with $\delta <0$ which means that 
the exponent $\frac{1}{q}-\frac{n-2}{2n}$ in (\ref{11.3}) is sharp.\abs
Let us mention here that for $p>1$ problem (\ref{0}) can be rewritten using
the substitution $v:=u^{1-p}$ as a Cauchy problem for the super-fast
diffusion equation given by
\be{0fast}
	\left\{ \begin{array}{l}
	v_t=\nabla \cdot (v^{m-1} \nabla v), \qquad x\in\R^n, \ t>0, \\[1mm]
	v(x,0)=v_0(x):=u_0^{1-p}, \qquad x\in\R^n,
 	\end{array} \right.
\ee
where $m=-\frac{1}{p-1}<0$ and $v_0\in C^0(\R^n)$ is such that $v_0(x)\to \infty$
as $|x|\to\infty$. Of course, our results on decay rates for (\ref{0}) can be rephrased
as results on growth rates of $\inf_{x\in\R^n} v(x,t)$ for (\ref{0fast}) in an evident manner.\abs
The equation $u_t=u\Delta u$, as corresponding to the borderline case $p=1$ in (\ref{0}),
occurs in the study of nonlinear transport
phenomena (\cite{CGR}), soil freezing processes (\cite{RAP}) and magma
solidification (\cite{ROCD}), for example.\abs
Gagliardo-Nirenberg inequalities were used before to obtain 
different results on asymptotic behavior of solutions of nonlinear
diffusion equations as in (\ref{0fast}) for various ranges of $m$, see
\cite{BBDGV, dPD, DT1, DT2}, for instance. For contexts where (\ref{0fast}) with
$m<0$ arises,
as well as for summaries of results on this problem, we refer to \cite{DK,
Vbook}.\abs
This paper is organized in such a way that 
Theorem~\ref{theo11} will be the objective of Section~2, whereas our study on the decay rates of
solutions to (\ref{0}) can be found in Section~3.
\mysection{A Gagliardo-Nirenberg-type inequality}
\subsection{Properties of functions satisfying (H)}
With two exceptions formed by Lemma~\ref{lem7} and Lemma~\ref{lem13} in which indeed (H) is directly referred to, 
throughout the sequel we will make use of (H) only through the elementary 
consequences of (H) stated in the following two lemmata.
The first of these, to be applied in Lemma \ref{lem3} but also again in the proof of Theorem \ref{theo11},
inter alia entails a property of essentially superalgebraic growth of the function
${\cal L}$ therein, for our later purposes formulated by including the derivative ${\cal L}'$.
\begin{lem}\label{lem4}
  Assume that $s_0\in (0,1)$ and that ${\cal L}\in C^0([0,s_0)) \cap C^1((0,s_0))$ is positive and nondecreasing on $(0,s_0)$ 
  and such that {\rm (H)} is valid.	
  Then
  \be{4.2}
	\frac{s{\cal L}'(s)}{{\cal L}(s)} \le \frac{a}{\ln \frac{1}{s}}
	\qquad \mbox{for all } s\in (0,s_0),
  \ee
  and in particular,
  \be{4.3}
	\frac{s{\cal L}'(s)}{{\cal L}(s)} \to 0
	\qquad \mbox{as } s\searrow 0.
  \ee
\end{lem}
\proof
  Given $s\in (0,s_0)$, we have $s<1$ and thus in particular $s-s^{1+\lambda}>0$ for all $\lambda>0$, whence we may apply
  e.g.~l'Hospital's rule to see that
  \be{4.4}
	\lim_{\lambda\searrow 0} \frac{-\lambda s\ln s}{s-s^{1+\lambda}}
	= \lim_{\lambda\searrow 0} \frac{s\ln s}{s^{1+\lambda} \ln s} =1.
  \ee
  On the other hand, (\ref{4.1}) implies that
  \bas
	{\cal L}(s)-{\cal L}(s^{1+\lambda})
	\le (1+a\lambda){\cal L}(s^{1+\lambda}) - {\cal L}(s^{1+\lambda})
	= a\lambda {\cal L}(s^{1+\lambda})
	\qquad \mbox{for all } \lambda\in (0,\lambda_0),
  \eas
  so that 
  \be{4.5}
	\limsup_{\lambda\searrow 0} \frac{{\cal L}(s)-{\cal L}(s^{1+\lambda})}{-\lambda s\ln s} 
	\le \limsup_{\lambda\searrow 0} \frac{a\lambda {\cal L}(s^{1+\lambda})}{-\lambda s \ln s}
	= \frac{a{\cal L}(s)}{s\ln \frac{1}{s}} \, ,
  \ee
  because ${\cal L}$ is continuous.
  As moreover ${\cal L}$ is even differentiable at $s$, combining (\ref{4.4}) with (\ref{4.5}) we thus infer that
  \bas
	{\cal L}'(s)
	&=& \lim_{\lambda\searrow 0} \frac{{\cal L}(s) - {\cal L}(s^{1+\lambda})}{s-s^{1+\lambda}} 
	= \lim_{\lambda\searrow 0} \bigg\{ \frac{{\cal L}(s)-{\cal L}(s^{1+\lambda})}{-\lambda s \ln s}
	\cdot \frac{-\lambda s\ln s}{s-s^{1+\lambda}} \bigg\} \\
	&\le& \bigg\{ \limsup_{\lambda\searrow 0} \frac{{\cal L}(s)-{\cal L}(s^{1+\lambda})}{-\lambda s \ln s} \bigg\}
	 \lim_{\lambda\searrow 0} \frac{-\lambda s \ln s}{s-s^{1+\lambda}} 
	= \frac{a{\cal L}(s)}{s\ln \frac{1}{s}} \, ,
  \eas
  which yields (\ref{4.2}) and thereby also entails (\ref{4.3}) due to the monotonicity of ${\cal L}$.
\qed
Apart from the latter, in Lemma \ref{lem3} we shall also make use of (H) through the following conclusion 
which is weaker than that in Lemma \ref{lem4} and actually satisfied also by any function ${\cal L}$ 
with precise algebraic behavior near the origin.
\begin{lem}\label{lem5}
  Let $s_0\in (0,1)$ and ${\cal L}\in C^0([0,s_0)) \cap C^1((0,s_0))$ be positive and nondecreasing on $(0,s_0)$ 
  and such that {\rm (H)} holds.	
  Then for any $d\in (0,1)$ there exists $C>0$ such that
  \be{5.2}
	{\cal L}(ds) \ge C{\cal L}(s)
	\qquad \mbox{for all } s\in (0,s_0).
  \ee
\end{lem}
\proof
  Writing $c_1:=a(\ln \frac{1}{s_0})^{-1}$, from (\ref{4.2}) we know that
  \bas
	\frac{{\cal L}'(s)}{{\cal L}(s)} \le \frac{c_1}{s}
	\qquad \mbox{for all } s\in (0,s_0),
  \eas
  which on integration shows that for fixed $d\in (0,1)$ and any $s\in (0,s_0)$ we have
  \bas
	\ln \frac{{\cal L}(s)}{{\cal L}(ds)} 
	\le \ln \frac{s^{c_1}}{(ds)^{c_1}},
  \eas
  so that thus (\ref{5.2}) holds with $C:=d^{c_1}$.
\qed
\subsection{Interpolation in Lebesgue spaces for rapidly decaying functions}
Now a major step toward our derivation of Theorem \ref{theo11} will be accomplished in the next lemma, 
the outcome of which already anticipates the structure of the desired inequality in (\ref{11.3}) but
yet exclusively contains Lebesgue norms of the considered function itself, rather than its gradient.
Accordingly, in the case of algebraic ${\cal L}$ given by ${\cal L}(s)=s^r$ with $r>0$,
the achieved estimate (\ref{3.2}) essentially reduces to a H\"older-type interpolation.
\begin{lem}\label{lem3}
  Assume that $s_0\in (0,1)$ and ${\cal L}\in C^0([0,\infty)) \cap C^1((0,s_0))$ is nonnegative, 
  nondecreasing and such that {\rm (H)} holds.
  Then for any choice of $n\ge 1$, $q_\star>0$, $q\in (0,q_\star)$ and $K>0$ one can find
$C=C(n,q,q_\star,K)>0$
  with the property that if $0\not\equiv \varphi \in L^{q_\star}(\R^n)$ is nonnegative and such that
  \be{3.1}
	\int_{\R^n} {\cal L}(\varphi) \le K,
  \ee
  then the inequality
  \be{3.2}
	\|\varphi\|_{L^q(\R^n)} \le C \|\varphi\|_{L^{q_\star}(\R^n)} 
	\bigg\{ {\cal L}^{-(\frac{1}{q}-\frac{1}{q_\star})} \Big(\|\varphi\|_{L^{q_\star}(\R^n)}^2 \Big) + 1 \bigg\}
  \ee
  holds.
\end{lem}
\proof
  We first recall the outcome of Lemma~\ref{lem4} 
  to find $s_1\in (0,1)$ such that $\frac{s{\cal L}'(s)}{{\cal L}(s)} \le \min\{q,\frac{q_\star}{2}\}$
  for all $s\in (0,s_1)$, which implies that
  \bea{3.5}
	\frac{d}{ds} \Big\{ s^{-q} {\cal L}(s)\Big\}
	= s^{-q-1} \Big\{ s{\cal L}'(s) - q{\cal L}(s) \Big\} 
	\le 0
	\qquad \mbox{for all } s\in (0,s_1),
  \eea
  and that similarly $\frac{d}{ds} \Big(s^{-\frac{q_\star}{2}} {\cal L}(s)\Big) \le 0$ for all $s\in (0,s_1)$.
  On integration, the latter entails that
  \bas
	s^{-\frac{q_\star}{2}} {\cal L}(s) \ge c_1:=s_1^{-\frac{q_\star}{2}} {\cal L}(s_1)
	\qquad \mbox{for all } s\in (0,s_1),
  \eas
  whence
  \be{3.6}
	{\cal L}(s) \ge c_1 s^\frac{q_\star}{2}
	\qquad \mbox{for all } s\in (0,s_1).
  \ee
  We now fix positive numbers $D$ and $d$ such that
  \be{3.7}
	D^{q_\star} \ge K
  \ee
  as well as
  \be{3.8}
	d<1
	\qquad \mbox{and} \qquad
	d \le c_1^{\frac{2}{q_\star}}D^{-2},
  \ee
  and thereafter apply Lemma~\ref{lem5} to choose $c_2>0$ satisfying
  \be{3.9}
	{\cal L}(ds) \ge c_2 {\cal L}(s)
	\qquad \mbox{for all } s\in (0,s_1).
  \ee
  We finally pick $s_2>0$ small enough such that $s_2 \le s_1$ and
  \be{3.10}
	s_2 \le \frac{1}{D} s_1 {\cal L}^\frac{1}{q_\star} (ds_1^2),
  \ee
  and suppose that $\varphi\in L^{q_\star}(\R^n)$ is nonnegative and such that $\varphi\not\equiv 0$ and 
  $\int_{\R^n} {\cal L}(\varphi) \le K$.
  Then 
  \be{3.11}
	B:=D^{q_\star-q} \|\varphi\|_{L^{q_\star}(\R^n)}^{-(q_\star-q)} 
	{\cal L}^{-\frac{q_\star-q}{q_\star}} \Big(d\|\varphi\|_{L^{q_\star}(\R^n)}^2 \Big)
  \ee
  is a well-defined positive number, and we first consider the case when
  \be{3.12}
	B^{-\frac{1}{q_\star-q}} \ge s_2,
  \ee
  in which we estimate the expression on the left-hand side of (\ref{3.2}) according to
  \be{3.13}
	\|\varphi\|_{L^q(\R^n)}
	\le \|\varphi\|_{L^q(\{\varphi\ge s_2\})}
	+ \|\varphi\|_{L^q(\{\varphi<s_2\})}.
  \ee
  Here, in view of (\ref{3.1}) and the monotonicity of ${\cal L}$ we see that
  \bas
	K \ge \int_{\R^n} {\cal L}(\varphi) \ge {\cal L}(s_2) \Big| \{\varphi \ge s_2\} \Big|,
  \eas
  and hence employing the H\"older inequality we obtain
  \bas
	\int_{\{\varphi\ge s_2\}} \varphi^q
	&\le& \bigg( \int_{\{\varphi\ge s_2\}} \varphi^{q_\star} \bigg)^\frac{q}{q_\star}
	\Big| \{\varphi \ge s_2\} \Big|^\frac{q_\star-q}{q_\star} 
	\le \bigg( \int_{\R^n} \varphi^{q_\star} \bigg)^\frac{q}{q_\star} 
	\bigg(\frac{K}{{\cal L}(s_2)} \bigg)^\frac{q_\star-q}{q_\star},
  \eas
  that is,
  \be{3.14}
	\|\varphi\|_{L^q(\{\varphi\ge s_2\})} \le c_3 \|\varphi\|_{L^{q_\star}(\R^n)}
  \ee
  holds with $c_3:=\Big(\frac{K}{{\cal L}(s_2)} \Big)^{\frac{1}{q}-\frac{1}{q_\star}}$.\\
  To control the second summand on the right of (\ref{3.13}) we first make use of the monotonicity property expressed
  in (\ref{3.5}) to see that since $s_2 \le s_1$ we have
  \bas
	\frac{\varphi^q(x)}{{\cal L}(\varphi(x))} \le c_4:=\frac{s_2^q}{{\cal L}(s_2)}
	\qquad \mbox{for all } x\in \{\varphi<s_2\}
  \eas
  and thus
  \bas
	\int_{\{\varphi<s_2\}} \varphi^q \le c_4 \int_{\{\varphi<s_2\}} {\cal L}(\varphi)
	\le c_4 K
  \eas
  again by (\ref{3.1}).
  In conjunction with (\ref{3.14}) and (\ref{3.13}), this shows that if (\ref{3.12}) is valid then
  \bas
	\|\varphi\|_{L^q(\R^n)}
	\le c_3\|\varphi\|_{L^{q_\star}(\R^n)} + c_4 K,
  \eas
  from which (\ref{3.2}) immediately follows in this case.\abs
  Conversely, if instead of (\ref{3.12}) we have
  \be{3.144}
	B^{-\frac{1}{q_\star-q}} < s_2,
  \ee
  then we first note that necessarily
  \be{3.15}
	\|\varphi\|_{L^{q_\star}(\R^n)} \le s_1,
  \ee
  because if this was false then by definition (\ref{3.11}) of $B$ and once more due to the monotonicity of ${\cal L}$ we would 
  obtain
  \bas
	s_2 > B^{-\frac{1}{q_\star-q}}
	= \frac{1}{D} \|\varphi\|_{L^{q_\star}(\R^n)} 
	{\cal L}^\frac{1}{q_\star} \Big( d\|\varphi\|_{L^{q_\star}(\R^n)}^2 \Big) 
	\ge \frac{1}{D} s_1 {\cal L}^\frac{1}{q_\star} (ds_1^2),
  \eas
  which is absurd in view of (\ref{3.10}).\\
  We shall next verify that for each $x\in\R^n$ we have
  \be{3.16}
	\varphi^q(x) \le B \varphi^{q_\star}(x)
	+{\cal B}  {\cal L}(\varphi(x)),\qquad
	{\cal B}:= \frac{1}{B^\frac{q}{q_\star-q} {\cal L}\left(B^{-\frac{1}{q_\star-q}}\right)} \, ,
  \ee
  which is obvious if $\varphi^p(x) \le B\varphi^{q_\star}(x)$, that is, if $\varphi(x) \ge B^{-\frac{1}{q_\star-q}}$.
  If $\varphi(x) < B^{-\frac{1}{q_\star-q}}$, however, then by our current assumption (\ref{3.144}) on $B$ we have
  $\varphi(x)<s_2 \le s_1$, and hence again the monotonicity property (\ref{3.5}) implies that
  \bas
	\frac{\varphi^q(x)}{{\cal L}(\varphi(x))}
	\le \frac{\left(B^{-\frac{1}{q_\star-q}}\right)^q}{{\cal L}\left(B^{-\frac{1}{q_\star-q}}\right)}
	={\cal B}
  \eas
  which completes the proof of (\ref{3.16}).\\
  Now integrating (\ref{3.16}) we find that
  \bea{3.17}
	\int_{\R^n} \varphi^q
	&\le& B \int_{\R^n} \varphi^{q_\star}
	+ {\cal B} K,
  \eea
  and we claim that our choice of $B$ ensures that herein
  \be{3.18}
{\cal B} K	\le B \int_{\R^n} \varphi^{q_\star}.
  \ee
  Indeed, by (\ref{3.11}) we have
  \bea{3.19}
	\frac{B \int_{\R^n} \varphi^{q_\star}}{{\cal B} K}
	&=& \frac{1}{K} B^\frac{q_\star}{q_\star-q} {\cal L}\Big( B^{-\frac{1}{q_\star-q}} \Big) 
	\int_{\R^n} \varphi^{q_\star} \nn\\
	&=& \frac{1}{K} \bigg\{ D^{q_\star-q} \|\varphi\|_{L^{q_\star}(\R^n)}^{-(q_\star-q)} 
	{\cal L}^{-\frac{q_\star-q}{q_\star}} \Big(d\|\varphi\|_{L^{q_\star}(\R^n)}^2 \Big) \bigg\}^\frac{q_\star}{q_\star-q} \times 
	\nn\\
	& & \hspace*{10mm}
	\times {\cal L} \Bigg(  \bigg\{ D^{q_\star-q} \|\varphi\|_{L^{q_\star}(\R^n)}^{-(q_\star-q)}
	{\cal L}^{-\frac{q_\star-q}{q_\star}} \Big(d\|\varphi\|_{L^{q_\star}(\R^n)}^2 \Big) \bigg\}^{-\frac{1}{q_\star-q}} \Bigg)
	 \int_{\R^n} \varphi^{q_\star} \nn\\
	&=& \frac{D^{q_\star}}{K}
	 \frac{{\cal L}\Big( \frac{1}{D} \|\varphi\|_{L^{q_\star}(\R^n)} 
	{\cal L}^\frac{1}{q_\star} \big( d\|\varphi\|_{L^{q_\star}(\R^n)}^2 \big) \Big)}
	{{\cal L}\big( d \|\varphi\|_{L^{q_\star}(\R^n)}^2 \big)},
  \eea
  where thanks to (\ref{3.8}) and (\ref{3.15}) we know that
  \be{3.20}
	d\|\varphi\|_{L^{q_\star}(\R^n)}^2 < \|\varphi\|_{L^{q_\star}(\R^n)}^2
	\le s_1^2 \le s_1,
  \ee
  so that in particular, by (\ref{3.6}) and (\ref{3.8}), 	
  \bas
	\frac{1}{D} \|\varphi\|_{L^{q_\star}(\R^n)}
	{\cal L}^\frac{1}{q_\star} \Big(d\|\varphi\|_{L^{q_\star}(\R^n)}^2 \Big)
	&\ge& \frac{1}{D} \|\varphi\|_{L^{q_\star}(\R^n)} 
	c_1^\frac{1}{q_\star} \Big(d\|\varphi\|_{L^{q_\star}(\R^n)}^2 \Big)^\frac{1}{2} \\
	&=& \frac{c_1^\frac{1}{q_\star} \sqrt{d}}{D} \|\varphi\|_{L^{q_\star}(\R^n)}^2 
	\ge d\|\varphi\|_{L^{q_\star}(\R^n)}^2.
  \eas
  Once more by means of the monotonicity of ${\cal L}$,
  from (\ref{3.19}) and (\ref{3.7}) we thus obtain that
  \bas
	\frac{B \int_{\R^n} \varphi^{q_\star}}{{\cal B} K}
	\ge \frac{D^{q_\star}}{K} 
	\frac{{\cal L}\big( d\|\varphi\|_{L^{q_\star}(\R^n)}^2 \big)}
		{{\cal L}\big(d\|\varphi\|_{L^{q_\star}(\R^n)}^2\big)}
	\ge 1.
  \eas
  Having thereby proved (\ref{3.18}), we may use this to infer from (\ref{3.17}) that according to our definition (\ref{3.11})
  of $B$,
  \bas
	\int_{\R^n} \varphi^q
	&\le& 2B \int_{\R^n} \varphi^{q_\star} 
	= 2  \bigg\{ D^{q_\star-q} \|\varphi\|_{L^{q_\star}(\R^n)}^{-(q_\star-q)} 
	{\cal L}^{-\frac{q_\star-q}{q_\star}} \Big(d\|\varphi\|_{L^{q_\star}(\R^n)}^2 \Big) \bigg\} 
	 \int_{\R^n} \varphi^{q_\star} \\
	&=& 2D^{q_\star-q} \|\varphi\|_{L^{q_\star}(\R^n)}^q
	{\cal L}^{-\frac{q_\star-q}{q_\star}} \Big( d\|\varphi\|_{L^{q_\star}(\R^n)}^2 \Big),
  \eas
  that is,
  \bas
	\|\varphi\|_{L^q(\R^n)}
	\le (2D^{q_\star-q})^\frac{1}{q} \|\varphi\|_{L^{q_\star}(\R^n)}
	{\cal L}^{-(\frac{1}{q}-\frac{1}{q_\star})} \Big( d\|\varphi\|_{L^{q_\star}(\R^n)}^2 \Big).
  \eas
  Again making use of (\ref{3.20}) in estimating
  \bas
	{\cal L}\Big(d\|\varphi\|_{L^{q_\star}(\R^n)}^2 \Big)
	\ge c_2 {\cal L} \Big( \|\varphi\|_{L^{q_\star}(\R^n)}^2 \Big)
  \eas
  by means of (\ref{3.9}), from this we readily derive (\ref{3.2}) also in the case when (\ref{3.144}) holds.
\qed
\subsection{Proof of Theorem~\ref{theo11}}
Combining Lemma \ref{lem3} with the Gagliardo-Nirenberg inequality in its well-known form, by once more
making use of Lemma \ref{lem4} we can now establish our main result on interpolation for rapidly 
decreasing functions.\abs
\proofc of Theorem~\ref{theo11}.\quad
  Given $0\not\equiv \varphi\in W^{1,2}(\R^n)$ such that (\ref{11.1}) holds, we first note that if
  $\|\varphi\|_{L^q(\R^n)} \le \|\nabla\varphi\|_{L^2(\R^n)}$, then abbreviating $\alpha:=\frac{1}{q}-\frac{n-2}{2n}$
  we can estimate
  \be{11.22}
	\frac{\|\varphi\|_{L^q(\R^n)}}{\|\nabla\varphi\|_{L^2(\R^n)}  {\cal L}^{-\alpha}(\|\nabla\varphi\|_{L^2(\R^n)})}
	\le {\cal L}_\infty^\alpha  \frac{\|\varphi\|_{L^q(\R^n)}}{\|\nabla\varphi\|_{L^2(\R^n)}} \le {\cal L}_\infty^\alpha
  \ee
  with ${\cal L}_\infty:=\|{\cal L}\|_{L^\infty((0,\infty))}$ being finite due to the boundedness of ${\cal L}$.\\
  We are thus left with the case when
  \be{11.23}
	\|\varphi\|_{L^q(\R^n)} > \|\nabla\varphi\|_{L^2(\R^n)},
  \ee
  in which using that $q\in (0,\frac{2n}{(n-2)_+})$ we can fix a number $q_\star\ge 1$ such that
  $q_\star>q$ and $q_\star\le \frac{2n}{(n-2)_+}$, so that an application of
  Lemma~\ref{lem3} yields $c_1>0$ fulfilling
  \be{4.39}
	\|\varphi\|_{L^q(\R^n)} \le c_1\|\varphi\|_{L^{q_\star}(\R^n)}  
	\bigg\{ {\cal L}^{-\gamma} \Big(\|\varphi\|_{L^{q_\star}(\R^n)}^2\Big) +1 \bigg\},
  \ee
  where $\gamma:=\frac{1}{q}-\frac{1}{q_\star}>0$.
  Here by means of the standard Gagliardo-Nirenberg inequality we can find $c_2\ge 1$ such that
  \be{4.40}
	\|\varphi\|_{L^{q_\star}(\R^n)}
	\le c_2\|\nabla\varphi\|_{L^2(\R^n)}^\theta \|\varphi\|_{L^q(\R^n)}^{1-\theta}
  \ee
  with $\theta:=\frac{\frac{n}{q}-\frac{n}{q_\star}}{1+\frac{n}{q}-\frac{n}{q_\star}}\in (0,1]$,
  and in order to make appropriate use of this on the right-hand side of (\ref{4.39}) we recall Lemma~\ref{lem4} 
  to pick $s_1>0$ satisfying
  \bas
	\frac{s{\cal L}'(s)}{{\cal L}(s)} \le \frac{1}{2\gamma}
	\qquad \mbox{for all } s\in (0,s_1),
  \eas
  which, namely, warrants that for
  \be{11.25}
	\rho(\sigma):=c_1 \sigma  \Big\{ {\cal L}^{-\gamma}(\sigma^2) +1\Big\}, \qquad \sigma>0,
  \ee
  we have
  \bas
	\rho'(\sigma)
	= c_1 + c_1 {\cal L}^{-\gamma-1}(\sigma^2)  \Big\{{\cal L}(\sigma^2)-2\gamma \sigma^2 {\cal L}'(\sigma^2)\Big\} \ge 0
	\qquad \mbox{for all } \sigma\in (0,\sqrt{s_1}).
  \eas
  Therefore, if $c_2\|\nabla\varphi\|_{L^2(\R^n)}^\theta \|\varphi\|_{L^q(\R^n)}^{1-\theta} < \sqrt{s_1}$
  then we obtain from (\ref{4.40}) and (\ref{4.39}) that
  \bas
	\|\varphi\|_{L^q(\R^n)}
	&\le& \rho \Big( c_2\|\nabla\varphi\|_{L^2(\R^n)}^\theta \|\varphi\|_{L^q(\R^n)}^{1-\theta}\Big) \\
	&\le& c_1 c_2 \|\nabla\varphi\|_{L^2(\R^n)}^\theta \|\varphi\|_{L^q(\R^n)}^{1-\theta} 
	\bigg\{ {\cal L}^{-\gamma} \Big( c_2^2 \|\nabla\varphi\|_{L^2(\R^n)}^{2\theta} \|\varphi\|_{L^q(\R^n)}^{2(1-\theta)}
		\Big) +1 \bigg\}
  \eas
  and hence
  \bea{11.24}
	\|\varphi\|_{L^q(\R^n)} 
	&\le& (c_1 c_2)^\frac{1}{\theta} \|\nabla\varphi\|_{L^2(\R^n)} 
	\bigg\{ {\cal L}^{-\gamma} \Big( c_2^2 \|\nabla\varphi\|_{L^2(\R^n)}^{2\theta} \|\varphi\|_{L^q(\R^n)}^{2(1-\theta)}
		\Big) +1 \bigg\}^\frac{1}{\theta} \nn\\
	&\le& (c_1 c_2)^\frac{1}{\theta} (1+{\cal L}_\infty^\gamma)^\frac{\gamma}{\theta} 
	\|\nabla\varphi\|_{L^2(\R^n)}		
	 {\cal L}^{-\frac{\gamma}{\theta}} \Big(c_2^2 \|\nabla\varphi\|_{L^2(\R^n)}^{2\theta} 
	\|\varphi\|_{L^q(\R^n)}^{2(1-\theta)} \Big) \nn\\
	&\le& (c_1 c_2)^\frac{1}{\theta} (1+{\cal L}_\infty^\gamma)^\frac{\gamma}{\theta} 
	\|\nabla\varphi\|_{L^2(\R^n)}		
	 {\cal L}^{-\frac{\gamma}{\theta}} \Big(\|\nabla\varphi\|_{L^2(\R^n)}^2 \Big),
  \eea
  because by definition of ${\cal L}_\infty$ we can estimate
  \bas
	1 \le {\cal L}_\infty^\gamma  {\cal L}^{-\gamma} \Big(c_2^2 \|\nabla\varphi\|_{L^2(\R^n)}^{2\theta} 
	\|\varphi\|_{L^q(\R^n)}^{2(1-\theta)} \Big),
  \eas
  and because (\ref{11.23}) along with our restriction $c_2\ge 1$ implies that
  \bas
	{\cal L}\Big(c_2^2 \|\nabla\varphi\|_{L^2(\R^n)}^{2\theta} 
	\|\varphi\|_{L^q(\R^n)}^{2(1-\theta)} \Big) 
	\ge {\cal L}\Big(\|\nabla\varphi\|_{L^2(\R^n)}^2 \Big).
  \eas
  If, conversely, $c_2\|\nabla\varphi\|_{L^2(\R^n)}^\theta \|\varphi\|_{L^q(\R^n)}^{1-\theta} \ge \sqrt{s_1}$,
  then (\ref{11.23}) entails that
  \be{11.26}
	\|\varphi\|_{L^q(\R^n)} \ge c_3:=\frac{\sqrt{s_1}}{c_2}.
  \ee
  In view of the fact that the function $\rho$ from (\ref{11.25}) satisfies $\rho(\sigma)\to 0$ as $\sigma\to 0$, 
  we can pick $\sigma_1>0$ such that $\rho(\sigma)<c_3$ for all $\sigma\in (0,\sigma_1)$, so that using 
  the inequality in (\ref{4.39}) in the form $\|\varphi\|_{L^q(\R^n)} \le \rho(\|\varphi\|_{L^{q_\star}(\R^n)})$,
  we infer from (\ref{11.26}) that necessarily $\|\varphi\|_{L^{q_\star}(\R^n)} \ge \sigma_1$.
  In conjunction with (\ref{4.39}), the monotonicity of ${\cal L}$ and (\ref{4.40}), however, this implies that writing
  $c_4:={\cal L}^{-\gamma}(\sigma_1^2)+1$ we have
  \bas
	\|\varphi\|_{L^q(\R^n)} 
	\le c_1 c_4 \|\varphi\|_{L^{q_\star}(\R^n)}
	\le c_1 c_2 c_4 \|\nabla\varphi\|_{L^2(\R^n)}^\theta \|\varphi\|_{L^q(\R^n)}^{1-\theta}
  \eas
  and thus 
  \bas
	\|\varphi\|_{L^q(\R^n)} \le (c_1 c_2 c_4)^\frac{1}{\theta} \|\nabla\varphi\|_{L^2(\R^n)},
  \eas
  whence proceeding as in (\ref{11.22}) we end up with the inequality
  \bas
	\frac{\|\varphi\|_{L^q(\R^n)}}{\|\nabla\varphi\|_{L^2(\R^n)}  {\cal L}^{-\alpha}(\|\nabla\varphi\|_{L^2(\R^n)})}
	\le (c_1 c_2 c_4)^\frac{1}{\theta} {\cal L}_\infty^\alpha
  \eas
  in this case. Together with (\ref{11.22}), (\ref{11.24}) and the observation that $\frac{\gamma}{\theta}=\alpha$,
  this establishes (\ref{11.3}).
\qed
\mysection{Decay estimates for solutions of $u_t=u^p \Delta u$}
\subsection{Preliminaries: Existence and approximation of solutions}
Next addressing the degenerate parabolic problem (\ref{0}) for $p\ge 1$,
in order to construct solutions thereof by approximation we follow \cite{fast_growth1} in considering
\be{0R}
	\left\{ \begin{array}{l}
	u_{Rt}=u_R^p \Delta u_R, \qquad x\in B_R, \ t>0, \\[1mm]
	u_R(x,t)=0, \qquad x\in\partial B_R, \ t>0, \\[1mm]
	u_R(x,0)=u_{0R}(x), \qquad x\in B_R,
 	\end{array} \right.
\ee
for $R>0$, where $u_{0R} \in C^3(\bar B_R)$ satisfies $0<u_{0R}<u_0$ in $B_R$ and $u_{0R}=0$ on $\partial B_R$ as well as
\be{conv_R}
	u_{0R} \nearrow u_0 \quad \mbox{in $\R^n$ \qquad as } R\nearrow \infty.
\ee
Moreover, for $\eps\in (0,1)$ we consider
\be{0Reps}
	\left\{ \begin{array}{l}
	u_{R\eps t}=u_{R \eps}^p \Delta u_{R\eps}, \qquad x\in B_R, \ t>0, \\[1mm]
	u_{R\eps}(x,t)=\eps, \qquad x\in\partial B_R, \ t>0, \\[1mm]
	u_{R\eps}(x,0)=u_{0R\eps}(x):=u_{0R}(x)+\eps, \qquad x\in B_R.
 	\end{array} \right.
\ee
Then the following basic statement has been shown in \cite{fast_growth1}.
\begin{lem}\label{lem202}
  Assume that $u_0\in C^0(\R^n) \cap L^\infty(\R^n)$ is positive.
  Then with $u_{0R}$ and $(u_{0\eps R})_{\eps\in (0,1)}$ as above, for each $\eps\in (0,1)$ the problem
  {\rm (\ref{0Reps})} possesses a global classical
  solution $u_{R\eps} \in C^0(\bar B_R \times [0,\infty)) \cap C^{2,1}(\bar B_R\times (0,\infty))$. 
  As $\eps\searrow 0$, we have $u_{R\eps}\searrow u_R$ with some positive classical solution
  $u_R\in C^0(\bar B_R \times [0,\infty)) \cap C^{2,1}(B_R\times (0,\infty))$ 
  of {\rm (\ref{0R})}. 
  Moreover, there exists a classical solution $u\in C^0(\R^n \times [0,\infty)) \cap C^{2,1}(\R^n\times (0,\infty))$
  of {\rm (\ref{0})} which is such that 
  \be{201.1}
	0<u(x,t)\le \|u_0\|_{L^\infty(\R^n)}	
	\qquad \mbox{ for all $x\in\R^n$ and } t\ge 0,
  \ee
  and that
  $u_R\nearrow u$ in $\R^n\times (0,\infty)$ as $R\nearrow \infty$.
  This solution is minimal in the sense that whenever $\tu\in C^0(\R^n \times [0,\infty)) \cap C^{2,1}(\R^n\times (0,\infty))$
  is positive and solves {\rm (\ref{0})} classically, we have $\tu\ge u$ in $\R^n\times [0,\infty)$.
\end{lem}
We note that in the special case when $u_0$ is radially symmetric around the origin and nonincreasing with respect
to $|x|$, we may and will assume that $u_{0R}$ has the same properties, which then, according to a standard argument involving
the comparison principle, are clearly inherited by $u_{R\eps}(\cdot,t)$ and hence also by $u_R(\cdot,t)$ for all
$t>0$.
\subsection{A Lyapunov functional ensuring persistence of fast spatial decay}\label{sect3.2}
To describe the large time asymptotics in (\ref{0}) using the above interpolation results,
let us first make sure that as a particular feature of the strong degeneracy in (\ref{0}) expressed
in our hypothesis $p\ge 1$, minimal solutions maintain the initial spatial decay.
Our general observation in this direction reads as follows.
%
%
\begin{lem}\label{lem9}
  Let $q>0$ and $s_0>0$, and suppose that the nonnegative and nondecreasing function 
  ${\cal L}\in C^0([0,s_0]) \cap C^2((0,s_0))$ has the property that
  \be{9.1}
	s{\cal L}''(s) \ge - \frac{3p+q-2}{p+q} \, {\cal L}'(s)
	\qquad \mbox{for all } s\in (0,s_0).
  \ee
  Then for any positive $u_0\in C^0(\R^n)$ satisfying $u_0^\frac{p+q}{2}<s_0$ in $\R^n$ and all $R>0$,
  there exists $\eps_0(R)\in (0,1)$ such that for each $\eps\in (0,\eps_0(R))$ the 
  solution $u_{R\eps}$ of {\rm (\ref{0Reps})} satisfies
  \be{9.2}
	\int_{B_R} {\cal L}\Big(u_{R\eps}^\frac{p+q}{2}(\cdot,t)\Big)
	\le \int_{B_R} {\cal L}\Big(u_{0\eps}^\frac{p+q}{2}\Big)
	\qquad \mbox{for all } t>0.
  \ee
\end{lem}
\proof
  Since $u_{0R}^\frac{p+q}{2} \le u_0^\frac{p+q}{2} < s_0$ in $\R^n$, for each $R>0$ we can find $\eps_0(R)\in (0,1)$
  such that $u_{0R\eps}^\frac{p+q}{2} < s_0$ in $\bar B_R$ for all $\eps\in (0,\eps_0(R))$.
  By comparison, this implies that the solution $u_{R\eps}$ of (\ref{0Reps}) satisfies
  $u_{R\eps}^\frac{p+q}{2} <0$ in $\bar B_R \times [0,\infty)$, so that (\ref{9.1}) applies to guarantee that
  \be{9.3}
	u_{R\eps}^\frac{p+q}{2} {\cal L}''\Big(u_{R\eps}^\frac{p+q}{2}\Big)
	\ge - \frac{3p+q-2}{p+q} \, {\cal L}'\Big(u_{R\eps}^\frac{p+q}{2}\Big)
	\qquad \mbox{in } \bar B_R \times [0,\infty).
  \ee
  Now from (\ref{0Reps}) we obtain that for all $R>0$ and $\eps\in (0,\eps_0(R))$,
  \bas
	\frac{2}{p+q} \frac{d}{dt} \int_{B_R} {\cal L}\Big( u_{R\eps}^\frac{p+q}{2}\Big)
	&=& \int_{B_R} u_{R\eps}^\frac{p+q-2}{2} {\cal L}'\Big(u_{R\eps}^\frac{p+q}{2}\Big) u_{R\eps t} 
	= \int_{B_R} u_{R\eps}^\frac{3p+q-2}{2} {\cal L}'\Big(u_{R\eps}^\frac{p+q}{2}\Big) \Delta u_{R\eps} \\
	&=& - \int_{B_R} \bigg\{
	\frac{p+q}{2} u_{R\eps}^\frac{4p+2q-4}{2} {\cal L}''\Big(u_{R\eps}^\frac{p+q}{2}\Big)
	+ \frac{3p+q-2}{2} u_{R\eps}^\frac{3p+q-4}{2} {\cal L}'\Big(u_{R\eps}^\frac{p+q}{2}\Big) \bigg\} 
	 |\nabla u_{R\eps}|^2 \\
	& & + \int_{\partial B_R} u_{R\eps}^\frac{3p+q-2}{2} {\cal L}'\Big(u_{R\eps}^\frac{p+q}{2}\Big) 
	\frac{\partial u_{R\eps}}{\partial\nu} \\[1mm]
	&\le& - \frac{p+q}{2} \int_{B_R} \bigg\{
	u_{R\eps}^\frac{p+q}{2} {\cal L}''\Big(u_{R\eps}^\frac{p+q}{2}\Big)
	+ \frac{3p+q-2}{p+q} {\cal L}'\Big(u_{R\eps}^\frac{p+q}{2}\Big) \bigg\}  
	u_{R\eps}^\frac{3p+q-4}{2} |\nabla u_{R\eps}|^2
  \eas
  for all $t>0$, 
  because ${\cal L}'\ge 0$ on $(0,s_0)$ and $\frac{\partial u_{R\eps}}{\partial\nu} \le 0$ on $\partial B_R\times (0,\infty)$
  due to the fact that $u_{R\eps} \ge \eps$ in $B_R \times (0,\infty)$ and $u_{R\eps}=\eps$
  on $\partial B_R \times (0,\infty)$.\\
  In view of (\ref{9.3}), however, this shows that $\frac{d}{dt} \int_{B_R} {\cal L}\Big(u_{R\eps}^\frac{p+q}{2}\Big) \le 0$
  for all $t>0$ and hence indeed
  \be{9.4}
	\int_{B_R} {\cal L}\Big(u_{R\eps}^\frac{p+q}{2}\Big)
	\le \int_{B_R} {\cal L}\Big( u_{0R\eps}^\frac{p+q}{2}\Big)
	\qquad \mbox{for all } t>0
  \ee
  whenever $R>0$ and $\eps\in (0,\eps_0(R))$.
%
\qed
When we choose ${\cal L}$ as a suitable power-type function, the above in particular
implies the control of the spatial $L^r$ quasi-norm
in the flavor
of (\ref{9.2}) for any $r\ge 1-p$, and hence for all positive $r$ whenever $p\ge 1$.
As the above reasoning shows, this conclusion actually extends to the not explicitly included cases $p=0$ and $p\in (0,1)$
corresponding to the heat equation and the porous medium equation, respectively, thus rediscovering
well-known Lyapunov-type properties of $\int_{\R^n} u^r$ for each $r\ge 1-p$ and any such $p$.
In view of our ambition to study solutions with fast spatial decay,
the essential role of our overall assumption $p\ge 1$ is underlined by the observation that the behavior
of these functionals drastically changes when $p<1$ and $r<1-p$.
Indeed, in the case $p=0$ it can directly be seen using explicit
solution representation through convolution with the Gauss kernel that for all nontrivial nonnegative initial
data in $L^1(\R^n)$ the corresponding functional $\int_{\R^n} u^r$ tends to $\infty$ as $t\to\infty$
for each $r\in (0,1)$; a similar conclusion can be drawn, e.g.~by using comparison from below with Barenblatt-type
self-similar solutions, when $p\in (0,1)$ and $r\in (0,1-p)$.
\abs
%
%
The requirement $p\ge 1$ guarantees that the above can actually be applied to 
functions ${\cal L}$ with a wide class of steepness properties near the origin.
Actually, instead of applying Lemma~\ref{lem9} directly, in our examples studied in
Corollaries~\ref{cor17}
and \ref{cor23} we will rather refer to the following weaker variant thereof.
\begin{lem}\label{lem99}
  Suppose that for some $s_0>0$,
  ${\cal L}\in C^0([0,s_0]) \cap C^2((0,s_0))$ is nonnegative and nondecreasing and such that
  \be{99.1}
	\frac{d}{ds} \Big(s{\cal L}'(s)\Big) \ge 0
	\qquad \mbox{for all } s\in (0,s_0).
  \ee
  Then for all $p\ge 1$ and $q>0$, the conclusion of Lemma~\ref{lem9} holds.
\end{lem}
\proof
  As (\ref{99.1}) implies that $s{\cal L}''(s)\ge - {\cal L}'(s)$ for all $s\in (0,s_0)$, 
  observing that this entails (\ref{9.1}) due to the fact that
  \bas
	\frac{3p+q-2}{p+q} = 1+\frac{2(p-1)}{p+q}\ge 1
	\qquad \mbox{for all $p\ge 1$ and } q>0,
  \eas
  we only need to apply Lemma~\ref{lem9}.
\qed
%
%
%
%
%
%
%
%
%
%
%
%
%
%
%
%
%
\subsection{Upper bounds in $L^q$ for $q>0$}
Having at hand the above information on conservation of spatial decay, we shall next address a statement
resembling that in Theorem \ref{theo14} but involving quasi-norms in $L^q(\R^n)$ for finite $q>0$.
Our result in this direction, to be achieved in Lemma \ref{lem12}, 
will be prepared by two lemmata, the first of them solves some transcendental inequalities involving
${\cal L}$ by once more explicitly referring to (H).
\begin{lem}\label{lem7}
  Assume that ${\cal L}\in C^0([0,\infty)) \cap C^1((0,\infty))$ is nondecreasing and nonnegative and satisfies
{\rm (H)},
  and let $\beta>\frac{1}{1+\lambda_0}$,
  $\gamma>0$ and $\delta_0>0$.
  Then there exists $C>0$ such that whenever $\eta>0$ is such that
  \be{7.1}
	\eta^\beta {\cal L}^\gamma(\eta) \le \delta
  \ee
  with some $\delta \in (0,\delta_0]$, then
  \be{7.2}
	\eta \le C \delta^\frac{1}{\beta} {\cal L}^{-\frac{\gamma}{\beta}}(\delta).
  \ee
\end{lem}
\proof
  Since $1+\lambda_0>\frac{1}{\beta}$, it is possible to find $\lambda\in (0,\lambda_0)$ such that still 
  $1+\lambda>\frac{1}{\beta}$, whence invoking (H) provides $s_1>0$ and $c_1>0$ such that
  \be{7.3}
	{\cal L}(s) \le c_1 {\cal L}(s^{1+\lambda})
	\qquad \mbox{for all } s\in (0,s_1).
  \ee
  We then pick $D>0$ large such that
  \be{7.4}
	D^\beta \ge c_1^\gamma
  \ee
  and
  \be{7.5}
	c_2 D \ge \delta_0^{1+\lambda-\frac{1}{\beta}},
  \ee
  where $c_2:={\cal L}^{-\frac{\gamma}{\beta}}(\delta_0)>0$.\\
  Now assuming (\ref{7.1}) to be valid for some $\eta>0$ and $\delta\in (0,\delta_0]$, we first consider the case when
  $\delta<s_1$, in which we claim that
  \be{7.6}
	\eta \le D\delta^\frac{1}{\beta} {\cal L}^{-\frac{\gamma}{\beta}}(\delta).
  \ee
  In fact, if on the contrary we had $\eta>D\delta^\frac{1}{\beta} {\cal L}^{-\frac{\gamma}{\beta}}(\delta)$, then
  by monotonicity of ${\cal L}$ we would have
  \bea{7.7}
	\frac{1}{\delta}\eta^\beta {\cal L}^\gamma(\eta)
	>\frac{1}{\delta}\Big(D\delta^\frac{1}{\beta} {\cal L}^{-\frac{\gamma}{\beta}}(\delta)\Big)^\beta 
		{\cal L}^\gamma \Big(D\delta^\frac{1}{\beta} {\cal L}^{-\frac{\gamma}{\beta}}(\delta)\Big)
	=\frac{ D^\beta}{{\cal L}^\gamma(\delta)}
	{\cal L}^\gamma\Big(D\delta^\frac{1}{\beta} {\cal L}^{-\frac{\gamma}{\beta}}(\delta)\Big).
  \eea
  Here since $\delta<s_1$ we may employ (\ref{7.3}) to estimate
  \be{7.8}
	{\cal L}^\gamma(\delta) \le c_1^\gamma(\delta^{1+\lambda}),
  \ee
  and again using the monotonicity of ${\cal L}$ we see that 
  ${\cal L}^{-\frac{\gamma}{\beta}}(\delta) \ge {\cal L}^{-\frac{\gamma}{\beta}}(\delta_0)=c_2$ and thus
  \be{7.9}
	{\cal L}^\gamma \Big( D\delta^\frac{1}{\beta} {\cal L}^{-\frac{\gamma}{\beta}}(\delta)\Big)
	\ge {\cal L}^\gamma \big(c_2 D\delta^\frac{1}{\beta}\big).
  \ee
  As our choice of $\lambda$ ensures that
  \bas
	\frac{\delta^{1+\lambda}}{c_2 D \delta^\frac{1}{\beta}}
	\le \frac{\delta_0^{1+\lambda-\frac{1}{\beta}}}{c_2 D} \le 1
  \eas
  by (\ref{7.5}), inserting (\ref{7.8}) and (\ref{7.9}) into (\ref{7.7}) and recalling (\ref{7.4}) therefore shows that
  \bas
	\frac{\eta^\beta {\cal L}^\gamma(\eta)}{\delta}
	> \frac{D^\beta}{c_1^\gamma} 
	\frac{{\cal L}^\gamma(c_2 D\delta^\frac{1}{\beta})}{{\cal L}^\gamma(c_2 D \delta^\frac{1}{\beta})} \ge 1.
  \eas
  This contradiction to (\ref{7.1}) warrants that indeed (\ref{7.6}) holds if $\delta<s_1$.\\
  However, if $s_1 \le \delta \le \delta_0$ then we observe that $\xi^\beta {\cal L}^\gamma(\xi)\to +\infty$ as $\xi\to + \infty$
  to verify that $\eta_0:=\sup\{\xi>0 \ | \ \xi^\beta {\cal L}^\gamma(\xi) \le \delta_0\}$ is well-defined and satisfies
  $\eta_0\ge \eta$ according to (\ref{7.1}).
  On the other hand, by definition of $c_2$ we have
  \bas
	\delta^{-\frac{1}{\beta}} {\cal L}^{-\frac{\gamma}{\beta}}(\delta)
	\ge \delta^\frac{1}{\beta} {\cal L}^{-\frac{\gamma}{\beta}}(\delta_0)
	= c_2 \delta^\frac{1}{\beta} \ge c_2 s_1^\frac{1}{\beta},
  \eas
  because $\delta \le \delta_ß$ and $\delta \ge s_1$.
  Consequently, in this case we obtain
  \bas
	\frac{\eta}{\delta^\frac{1}{\beta} {\cal L}^{-\frac{\gamma}{\beta}}(\delta)}
	\le \frac{\eta_0}{c_2 s_1^\frac{1}{\beta}},
  \eas
  and thus we all in all conclude that (\ref{7.2}) is valid if we let
  $C:=\max\Big\{D \, , \, \frac{\eta_0}{c_2 s_1^{1/\beta}} \Big\}$.
\qed
Another consequence of (H) used in Lemma \ref{lem12} states that for fixed nonnegative measurable and bounded 
$\varphi$, the family $({\cal L}(\varphi^r))_{r>0}$ either entirely belongs to $L^1(\R^n)$ or lies completely outside,
which clearly again reflects a strongly superalgebraic growth of ${\cal L}(s)$ near $s=0$.
\begin{lem}\label{lem13}
  Suppose that ${\cal L}\in C^0([0,\infty))$ is nonnegative and nondecreasing and such that {\rm (H)} is valid.
  Then for any nonnegative $\varphi\in L^\infty(\R^n)$ satisfying
  \be{13.1}
	\int_{\R^n} {\cal L}(\varphi)<\infty,
  \ee
  we have
  \be{13.2}
	\int_{\R^n} {\cal L}(\varphi^r) < \infty
	\qquad \mbox{for all } r>0.
  \ee
\end{lem}
\proof
  Without loss of generality assuming that $s_0\le 1$, we first note that since
  \bas
	\int_{\R^n} {\cal L}(\varphi) \ge \int_{\{\varphi\ge s_0^\frac{1}{r}\}} {\cal L}(\varphi)
	\ge {\cal L}(s_o^\frac{1}{r})  \Big| \{\varphi \ge s_0^\frac{1}{r}\}\Big|
  \eas
  by monotonicity of ${\cal L}$, (\ref{13.1}) asserts that $c_1:=|\{\varphi\ge s_0^\frac{1}{r}\}|$ is finite.
  Now in the case $r<1$ it is easy to see that there exist $k\in\N$ and $\lambda\in (0,\lambda_0)$ such that
  $r(1+\lambda)^k=1$, whence $k$ applications of (H) yield
  \bas
	{\cal L}(s^r) \le (1+a\lambda)^k {\cal L}\Big(s^{r(1+\lambda)^k}\Big) = (1+a\lambda)^k {\cal L}(s)
	\qquad \mbox{for all } s\in [0,s_0^\frac{1}{r}),
  \eas
  because for any such $s$ and each $j\in\{0,...,k-1\}$ we have $s^{r(1+\lambda)^j} < s_0^{(1+\lambda)^j} \le s_0$
  due to the fact that $s_0\le 1$. 
  Accordingly, again by monotonicity of ${\cal L}$ we obtain
  \bas
	\int_{\R^n} {\cal L}(\varphi^r)
	&=& \int_{\{\varphi \ge s_0^\frac{1}{r}\}} {\cal L}(\varphi^r)
	+ \int_{\{\varphi < s_0^\frac{1}{r}\}} {\cal L}(\varphi^r) \\
	&\le& {\cal L}\Big(\|\varphi\|_{L^\infty(\R^n)}^r\Big)  \Big| \{\varphi \ge s_0^\frac{1}{r}\} \Big|
	+ \int_{\{\varphi< s_0^\frac{1}{r}\}} (1+a\lambda)^k {\cal L}(\varphi) \\
	&\le& c_1 {\cal L}\Big(\|\varphi\|_{L^\infty(\R^n)}^r\Big) 
	+ (1+a\lambda)^k \int_{\R^n} {\cal L}(\varphi),
  \eas
  so that (\ref{13.1}) indeed implies (\ref{13.2}) in this case.\\
  If $r\ge 1$, however, we similarly estimate
  \bas
	\int_{\R^n} {\cal L}(\varphi^r)
	&=& \int_{\{\varphi\ge 1\}} {\cal L}(\varphi^r)
	+\int_{\{\varphi< 1\}} {\cal L}(\varphi^r) \\
	&\le& {\cal L}\Big(\|\varphi\|_{L^\infty(\R^n)}^r\Big)  \Big| \{\varphi\ge 1\} \Big|
	+\int_{\{\varphi< 1\}} {\cal L}(\varphi^r),
  \eas
  where clearly $s_0\le 1$ entails that $|\{\varphi \ge 1\}| \le c_1$, and where $r\ge 1$ implies that
  \bas
	\int_{\{\varphi< 1\}} {\cal L}(\varphi^r)
	\le \int_{\{\varphi<1\}} {\cal L}(\varphi)
	\le \int_{\R^n} {\cal L}(\varphi),
  \eas
  whence again (\ref{13.2}) results from (\ref{13.1}).
\qed
Using Theorem~\ref{theo11} along with Lemma~\ref{lem9},
we can now achieve an essential step toward Theorem~\ref{theo14} by
deriving a corresponding $L^q$ counterpart for solutions to the approximate system (\ref{0R}).
Indeed, our argument will be based on a refined examination of the
time evolution of $L^q$ quasi-norms along trajectories of (\ref{0R}), where unlike
in Lemma~\ref{lem9} we shall here rely on the interpolation property from Theorem~\ref{theo11} in gaining
a nontrivial estimate from below for the corresponding dissipation rate (cf.~(\ref{12.33}) and (\ref{12.7})).
\begin{lem}\label{lem12}
%
  Suppose that $s_0>0$ and that ${\cal L}\in C^0([0,\infty)) \cap C^2((0,s_0))$ is positive and nondecreasing on $(0,\infty)$
  such that {\rm (H)} is satisfied, and such that
  \be{12.1}
	s{\cal L}''(s) \ge - \frac{3p+q_0-2}{p+q_0} {\cal L}'(s)
	\qquad \mbox{for all } s\in (0,s_0)
  \ee
  with a certain $q_0>0$.
  Moreover, assume that $n\ge 3$ and that $u_0\in C^0(\R^n)$ is positive and such that
  \be{12.01}
	u_0<\min \Big\{ s_0^\frac{2}{p} \, , \, s_0^\frac{2}{p+q_0} \Big\}
	\qquad \mbox{in } \R^n
  \ee
  as well as
  \be{12.2}
	\int_{\R^n} {\cal L}(u_0) < \infty.
  \ee
  Then there exist $q_1 \in (0,q_0)$ with the property that for all $q\in (0,q_1)$ one can find 
  $t_0=t_0(q)>0$ and $C=C(q)>0$ such that the solution $u_R$ of (\ref{0R}) satisfies
  \be{12.3}
	\|u_R(\cdot,t)\|_{L^q(B_R)}
	\le C t^{-\frac{1}{p}} {\cal L}^{-\frac{np+2q}{npq}} \Big(\frac{1}{t}\Big)
	\qquad \mbox{for all } t\ge t_0.
  \ee
\end{lem}
\proof
  We fix any $q_1\in (0,q_0)$ such that
  \be{12.4}
	\frac{2q}{p+q} \le 1
	\quad \mbox{for all } q\in (0,q_1)
	\qquad \mbox{and} \qquad
	q_1 < p\lambda_0,
  \ee
  and given $q\in (0,q_1)$ we may combine (\ref{12.2}) with the outcome of Lemma~\ref{lem13} to see that
  \bas
	\int_{\R^n} {\cal L}\Big(u_0^\frac{p+q}{2}\Big)<\infty.
  \eas
  As (\ref{12.01}) implies that moreover $u_0^\frac{p+q}{2}<s_0$ in $\R^n$, from Lemma~\ref{lem9}
  we thus infer that for any $R>0$ one can find $\eps_0(R)\in (0,1)$ such that whenever $\eps\in (0,\eps_0(R))$,
  the solution of (\ref{0Reps}) satisfies
  \be{12.5}
	\int_{B_R} {\cal L}\Big( u_{R\eps}^\frac{p+q}{2}(\cdot,t)\Big) 
	\le \int_{B_R} {\cal L}\Big( u_{0R\eps}^\frac{p+q}{2}\Big)
	\qquad \mbox{for all } t>0.
  \ee
  Since the monotonicity of ${\cal L}$ ensures that
  \bas
	\int_{B_R} {\cal L}\Big(u_{0R\eps}^\frac{p+q}{2}\Big) \searrow 
	\int_{B_R} {\cal L}\Big( u_{0R}^\frac{p+q}{2}\Big)
	\qquad \mbox{as } \eps\searrow 0
  \eas
  by Beppo Levi's theorem,
  and that furthermore 
  \bas
	\int_{B_R} {\cal L}\Big( u_{0R}^\frac{p+q}{2}\Big)
	\le c_1:=\int_{\R^n} {\cal L}\Big( u_0^\frac{p+q}{2}\Big)
	\qquad \mbox{for all } R>0,
  \eas
  the inequality (\ref{12.5}) entails that for each $R>0$ we can fix $\eps_1(R)\in (0,\eps_0(R))$ such that for any
  $\eps\in (0,\eps_1(R))$ we have
  \be{12.6}
	\int_{B_R} {\cal L}\Big(u_{R\eps}^\frac{p+q}{2}(\cdot,t)\Big)	
	\le 2c_1
	\qquad \mbox{for all } t>0.
  \ee
  Now testing (\ref{0Reps}) by the smooth function $u_{R\eps}^{q-1}$ yields
  \bea{12.33}
	\frac{1}{q} \frac{d}{dt} \int_{B_R} u_{R\eps}^q
	&=& \int_{B_R} u_{R\eps}^{p+q-1} \Delta u_{R\eps} \nn\\
	&=& -(p+q-1) \int_{B_R} u_{R\eps}^{p+q-2} |\nabla u_{R\eps}|^2
	+ \int_{\partial B_R} u_{R\eps}^{p+q-1} \frac{\partial u_{R\eps}}{\partial\nu} \nn\\
	&\le& - \frac{4(p+q-1)}{(p+q)^2} \int_{B_R} \Big|\nabla u_{R\eps}^\frac{p+q}{2} \Big|^2
	\qquad \mbox{for all } t>0,
  \eea
  once again because $\frac{\partial u_{R\eps}}{\partial\nu} \le 0$ on $\partial B_R\times (0,\infty)$.
  In order to estimate the right-hand side herein by means of Theorem~\ref{theo11}, we observe that for each fixed $t>0$,
  the function $u_{R\eps}^\frac{p+q}{2}(\cdot,t)-\eps^\frac{p+q}{2} \in C^1(\bar B_R)$ is positive in $B_R$ and vanishes
  on $\partial B_R$, so that its trivial extension to all of $\R^n$ belongs to $W^{1,2}(\R^n)$.
  As
  \bas
	\int_{B_R} {\cal L}\Big(u_{R\eps}^\frac{p+q}{2}(\cdot,t)-\eps^\frac{p+q}{2}\Big)
	\le \int_{B_R} {\cal L}\Big( u_{R\eps}^\frac{p+q}{2}(\cdot,t)\Big)
	\le 2c_1
	\qquad \mbox{for all $t>0$ and each } \eps\in (0,\eps_1(R)),
  \eas
  Theorem~\ref{theo11} therefore becomes applicable so as to yield $c_2>0$ such that
  \bea{12.7}
	& & \hspace*{-20mm}
	\Big\| u_{R\eps}^\frac{p+q}{2}-\eps^\frac{p+q}{2}\Big\|_{L^\frac{2q}{p+q}(B_R)}^\frac{2q}{p+q} \nn\\
	&\le& c_2  \Bigg\{ \Big\| \nabla (u_{R\eps}^\frac{p+q}{2}-\eps^\frac{p+q}{2}) \Big\|_{L^2(B_R)}
	 {\cal L}^{-(\frac{p+q}{2q}-\frac{n-2}{2n})} \bigg( 
	\Big\| \nabla (u_{R\eps}^\frac{p+q}{2}-\eps^\frac{p+q}{2}) \Big\|_{L^2(B_R)}^2 \bigg) \Bigg\}^\frac{2q}{p+q} \nn\\
	&=& c_2  \bigg( \Big\| \nabla u_{R\eps}^\frac{p+q}{2}\Big\|_{L^2(B_R)}^2 \bigg)^\frac{q}{p+q}
	 {\cal L}^{-\frac{np+2q}{n(p+q)}} \bigg( \Big\| \nabla u_{R\eps}^\frac{p+q}{2} \Big\|_{L^2(B_R)}^2 \bigg)
	\qquad \mbox{for all } t>0,
  \eea
  where on the left-hand side we may use that thanks to the first restriction in (\ref{12.4}) we have
  $(x+y)^\frac{2q}{p+q} \le x^\frac{2q}{p+q}+y^\frac{2q}{p+q}$ for all $x\ge 0$ and $y\ge 0$, so that
  \bea{12.8}
	\int_{B_R} u_{R\eps}^q
	&=& \int_{B_R} \bigg\{ \Big(u_{R\eps}^\frac{p+q}{2}-\eps^\frac{p+q}{2}\Big) + \eps^\frac{p+q}{2} 
		\bigg\}^\frac{2q}{p+q} \nn\\
	&\le& \Big\| u_{R\eps}^\frac{p+q}{2} - \eps^\frac{p+q}{2} \Big\|_{L^\frac{2q}{p+q}(B_R)}^\frac{2q}{p+q}
	+ |B_R|  \eps^q
	\qquad \mbox{for all } t>0.
  \eea
  Now to solve (\ref{12.7}) with respect to $\|\nabla u_{R\eps}^\frac{p+q}{2}\|_{L^2(B_R)}^2$, we abbreviate
  $\beta:=\frac{q}{p+q}$ and $\gamma:=\frac{np+2q}{n(p+q)}$ and note that since
  \bas
	\frac{d}{ds} \Big\{ s^\beta {\cal L}^{-\gamma}(s)\Big\}
	= s^{\beta-1} {\cal L}^{-\gamma-1}(s)  \Big\{ \beta {\cal L}(s)-\gamma s{\cal L}'(s)\Big\}
	\qquad \mbox{for all } s\in (0,s_0),
  \eas
  invoking Lemma~\ref{lem4} provides $s_1\in (0,s_0)$ such that
  \bas
	\frac{d}{ds} \Big\{ s^\beta {\cal L}^{-\gamma}(s)\Big\}
	>0
	\qquad \mbox{for all } s\in (0,s_1].
  \eas
  The function $\psi$ defined on $[0,\infty)$ by letting
  \bas
	\psi(s):=c_2 s^\beta \tl^{-\gamma}(s), \qquad s\ge 0,
  \eas
  with
  \bas
	\tl(s):=\left\{ \begin{array}{ll}
	{\cal L}(s), \qquad & s\in [0,s_1], \\[1mm]
	{\cal L}(s_1), & s>s_1,
	\end{array} \right.
  \eas
  therefore has the properties that $\psi'>0$ on $(0,\infty)\setminus \{s_1\}$ and $\psi(0)=0$ as well as
  $\psi(s)\to + \infty$ as $s\to\infty$, and since ${\cal L}$ is nondecreasing we moreover have
  $\psi(s)\ge c_2 s^\beta {\cal L}^{-\gamma}(s)$ for all $s\ge 0$.
  Accordingly, combining (\ref{12.7}) with (\ref{12.8}) shows that for all $R>0$ and $\eps\in (0,\eps_1(R))$,
  \bas
	y_{R\eps}(t):=\int_{B_R} u_{R\eps}^q(\cdot,t) - |B_R| \eps^q,
	\qquad t\ge 0,
  \eas
  satisfies
  \bas
	y_{R\eps}(t)
	&\le& c_2  \bigg(\Big\| \nabla u_{R\eps}^\frac{p+q}{2}(\cdot,t) \Big\|_{L^2(B_R)}^2 \bigg)^\beta
	 {\cal L}^{-\gamma} \bigg(\Big\| \nabla u_{R\eps}^\frac{p+q}{2}(\cdot,t) \Big\|_{L^2(B_R)}^2 \bigg) \\
	&\le& \psi \bigg(\Big\| \nabla u_{R\eps}^\frac{p+q}{2}(\cdot,t) \Big\|_{L^2(B_R)}^2 \bigg)
	\qquad \mbox{for all } t>0,
  \eas
  so that since $u_{R\eps}>\eps$ in $B_R\times (0,\infty)$ entails that $y_{R\eps}$ is positive, we may invert this 
  relation so as to achieve that
  \bas
	\int_{B_R} \Big|\nabla u_{R\eps}^\frac{p+q}{2}(\cdot,t)\Big|^2
	\ge \psi^{-1} \big(y_{R\eps}(t)\big)
	\qquad \mbox{for all } t>0,
  \eas
  where $\psi^{-1}$ denotes the inverse of $\psi$.
  Abbreviating $c_3:=\frac{4q(p+q-1)}{(p+q)^2}$, from (\ref{12.33}) we thus obtain the autonomous ODI
  \bas
	y_{R\eps}'(t) \le -c_3 \psi^{-1}\big(y_{R\eps}(t)\big)
	\qquad \mbox{for all } t>0,
  \eas
  which again by positivity of $y_{R\eps}$ can be integrated to see that
  \bas
	\int_{y_{R\eps}(0)}^{y_{R\eps}(t)} \frac{dy}{\psi^{-1}(y)} \le -c_3 t
	\qquad \mbox{for all } t>0,
  \eas
  whence by substituting $s:=\psi^{-1}(y)$ we obtain that
  \be{12.9}
	c_3 t
	\le \int_{y_{R\eps}(t)}^{y_{R\eps}(0)} \frac{dy}{\psi^{-1}(y)}
	= \int_{\psi^{-1}(y_{R\eps}(t))}^{\psi^{-1}(y_{R\eps}(0))} \frac{\psi'(s)}{s} ds
	\qquad \mbox{for all } t>0.
  \ee
  Since herein the monotone convergence $u_{R\eps} \searrow u_R$ warrants that for all $t\ge 0$ we have
  \bas
	y_{R\eps}(t) \to y_R(t):=\int_{B_R} u_R^q(\cdot,t)
	\qquad \mbox{as } \eps\searrow 0,
  \eas
  by continuity of $\psi^{-1}$ we infer on taking $(0,\eps_1(R))\ni \eps \searrow 0$ in (\ref{12.9}) that
  \be{12.10}
	c_3 t \le \int_{\psi^{-1}(y_R(t))}^{\psi^{-1}(y_R(0))} \frac{\psi'(s)}{s} ds
	\qquad \mbox{for all } t>0.
  \ee
  Here thanks to the monotonicity of $\tl$, for any $\os>0$ we can estimate
  \bas
	\psi'(s)
	&=& \beta c_2 s^{\beta-1} \tl^{-\gamma}(s) - \gamma c_2 s^\beta \tl^{-\gamma-1}(s) \tl'(s) 
	\le\beta c_2 s^{\beta-1} \tl^{-\gamma}(s) \\
	&\le& \beta c_2 \tl^{-\gamma}(\os)  s^{\beta-1}
	\qquad \mbox{for all } s\in (\os,\infty)\setminus \{s_1\},
  \eas
  so that (\ref{12.10}) along with the fact that $\beta=\frac{q}{p+q}<1$ implies that
  \bas
	c_3 t
	&\le& \beta c_2 \tl^{-\gamma}\Big(\psi^{-1}(y_R(t))\Big)  
		\int_{\psi^{-1}(y_R(t))}^{\psi^{-1}(y_R(0))} s^{\beta-2} ds 
	\le \beta c_2 \tl^{-\gamma}\Big(\psi^{-1}(y_R(t))\Big)  
		\int_{\psi^{-1}(y_R(t))}^\infty s^{\beta-2} ds \\
	&=& \frac{\beta c_2}{1-\beta} \tl^{-\gamma}\Big(\psi^{-1}(y_R(t))\Big)  
	\Big( \psi^{-1}(y_R(t))\Big)^{\beta-1}
	\qquad \mbox{for all } t>0,
  \eas
  that is, we have
  \bas
	\Big( \psi^{-1}(y_R(t))\Big)^{1-\beta}  \tl^{\gamma} \Big(\psi^{-1}(y_R(t))\Big)
	\le \frac{\beta c_2}{1-\beta}  \frac{1}{c_3 t}
	\le \frac{c_4}{t}
	\qquad \mbox{for all } t>0
  \eas
  with $c_4:=\max\{ \frac{\beta c_2}{(1-\beta)c_3} \, , \, 1\}$.\\
  Invoking Lemma~\ref{lem7}, we thus infer the existence of $t_1>0$ and $c_5\ge 1$ satisfying
  \bas
	\psi^{-1}(y_R(t))
	\le c_5  \Big(\frac{c_4}{t}\Big)^\frac{1}{1-\beta}
 \tl^{-\frac{\gamma}{1-\beta}} \Big(\frac{c_4}{t}\Big) 
	\le c_4^\frac{1}{1-\beta} c_5  t^{-\frac{1}{1-\beta}}
		\tl^{-\frac{\gamma}{1-\beta}} \Big(\frac{1}{t}\Big)
	\qquad \mbox{for all } t\ge t_1,
  \eas
  because $c_4\ge 1$ and $\tl$ is nondecreasing.\\
  Writing $c_6:=c_4^\frac{1}{1-\beta} c_5$, upon inversion we thereby obtain that
  \bea{12.11}
	y_R(t)
	&\le& \psi \bigg( c_6 t^{-\frac{1}{1-\beta}} \tl^{-\frac{\gamma}{1-\beta}}\Big(\frac{1}{t}\Big) \bigg) \nn\\
	&=& c_2  \bigg\{ c_6 t^{-\frac{1}{1-\beta}} \tl^{-\frac{\gamma}{1-\beta}}\Big(\frac{1}{t}\Big) \bigg\}^\beta
	 \tl^{-\gamma} \bigg( c_6 t^{-\frac{1}{1-\beta}} \tl^{-\frac{\gamma}{1-\beta}}\Big(\frac{1}{t}\Big) \bigg) \nn\\
	&=& c_2 c_6^\beta t^{-\frac{\beta}{1-\beta}} \tl^{-\frac{\beta\gamma}{1-\beta}} \Big(\frac{1}{t}\Big)
	 \tl^{-\gamma} \bigg( c_6 t^{-\frac{1}{1-\beta}} \tl^{-\frac{\gamma}{1-\beta}}\Big(\frac{1}{t}\Big) \bigg)
	\qquad \mbox{for all } t\ge t_1.
  \eea
  Here the last factor can be estimated for large $t$ by choosing $t_0\ge t_1$ large enough fulfilling
  \be{12.12}
	\frac{1}{t_0} \le s_1
	\qquad \mbox{and} \qquad
	t_0^{-\frac{1}{1-\beta}} \le s_1
  \ee
  as well as
  \be{12.13}
	c_6 \tl^{-\frac{\gamma}{1-\beta}} \Big(\frac{1}{t}\Big) \ge 1
	\qquad \mbox{for all } t\ge t_0,
  \ee
  where the latter is possible since $\tl$ is continuous with $\tl(0)=0$.
  Using (\ref{12.13}) and the second restriction in (\ref{12.12}) we thus infer that
  \be{12.14}
	\tl^{-\gamma} \bigg( c_6 t^{-\frac{1}{1-\beta}} \tl^{-\frac{\gamma}{1-\beta}}\Big(\frac{1}{t}\Big) \bigg)
	\le \tl^{-\gamma} \Big(t^{-\frac{1}{1-\beta}}\Big)
	= {\cal L}^{-\gamma} \Big( t^{-\frac{1}{1-\beta}}\Big)
	\qquad \mbox{for all } t\ge t_0,
  \ee
  where since $\frac{1}{1-\beta}=\frac{p+q}{p}<1+\lambda_0$ by (\ref{12.4}), and since
  $\frac{1}{t_0} \le s_1 < s_0$ by (\ref{12.12}), we may invoke (\ref{4.1}) to find that
  \bas
	{\cal L}\Big(\frac{1}{t}\Big)
	\le c_7 {\cal L}\Big(t^{-\frac{1}{1-\beta}}\Big)
	\qquad \mbox{for all } t\ge t_0
  \eas
  with $c_7:=1+\frac{aq}{p}$.
  Combining this with (\ref{12.14}), (\ref{12.11}) and again (\ref{12.10}), we therefore conclude that
  \bea{12.15}
	y_R(t)
	&\le& c_2 c_6^\beta t^{-\frac{\beta}{1-\beta}} {\cal L}^{-\frac{\beta\gamma}{1-\beta}} \Big(\frac{1}{t}\Big)
	 c_7^\gamma {\cal L}^{-\gamma}\Big(\frac{1}{t}\Big) \nn\\
	&=& c_2 c_6^\beta c_7^\gamma t^{-\frac{\beta}{1-\beta}} {\cal L}^{-\frac{\gamma}{1-\beta}}\Big(\frac{1}{t}\Big)
	\qquad \mbox{for all } t\ge t_0.
  \eea
  As $\frac{\beta}{1-\beta}=\frac{q}{p}$ and $\frac{\gamma}{1-\beta}=\frac{np+2q}{np}$, taking $q$-th roots on both
  sides of (\ref{12.15}) readily yields (\ref{12.3}).
\qed
\subsection{Upper bounds in $L^\infty$. Proof of Theorem~\ref{theo14}}
In order to prepare our deduction of spatially uniform estimates from the above inequalities involving $L^q$ seminorms, 
we recall the following well-known semi-convexity property (\cite{aronson}, \cite{fast_growth1}).
\begin{lem}\label{lem87}
  Let $R>0$. Then the solution of (\ref{0R}) from Lemma \ref{lem202} satisfies
  \bas
	\frac{u_{Rt}(x,t)}{u_R(x,t)} \ge -\frac{1}{pt}
	\qquad \mbox{for all } x\in B_R \mbox{ and } t>0.
  \eas
\end{lem}
For a radially symmetric and radially nondecreasing solution, namely, this entails conrollability of its 
spatial $L^\infty$ norm by its $L^q$ seminorm for arbitrarily small $q>0$.
\begin{lem}\label{lem89}
%
  Assume that $u_0\in C^0(\R^n)$ is positive, radially symmetric and nondecreasing
  with respect to $|x|$. Then for any $q>0$, the solution of (\ref{0R}) from Lemma \ref{lem202} satisfies
  \be{89.1}
	\|u_R(\cdot,t)\|_{L^\infty(B_R)} 
	\le \Bigg( \frac{2^{q+\frac{n(p-1)}{2}}n}{p^\frac{n}{2} \omega_n}\Bigg)^\frac{2}{np+2q}
	t^{-\frac{n}{np+2q}} \|u_R(\cdot,t)\|_{L^q(B_R)}^\frac{2q}{np+2q}
	\qquad \mbox{for all } t>0,
  \ee
  where $\omega_n:=n|B_1|$.
\end{lem}
\proof
  Without danger of confusion we may write $u(r,t)$ for $r=|x|\ge 0$, and given $t>0$ we then let
  \be{89.2}
	r_0\equiv r_0(t):=\sup \Big\{ r\in (0,R) \ \Big| \ u_R(r,t) \ge \frac{1}{2} u_R(0,t) \Big\},
  \ee
  noting that $r_0$ is well-defined due to the fact that $u_R(0,t)>u_R(R,t)=0$. 
  Now from Lemma \ref{lem87} and (\ref{0R}) we know that
  \bas
	u_R^{p-1} \Delta u_R = \frac{u_{Rt}}{u_R} \ge -\frac{1}{pt}
	\qquad \mbox{in } B_R
  \eas
  and hence, as $u_R(\cdot,t)$ clearly inherits the symmetry and monotonicity properties of $u_{0R}$ by the 
  maximum principle, 
  \bas
	\partial_r^2 u_R(r,t)
	&\ge& \partial_r^2 u_R(r,t) + \frac{n-1}{r} \partial_r u_R(r,t) 
	\ge - \frac{1}{pt} u_R^{1-p}(r,t) \\
	&\ge& - \frac{1}{pt}  \Big(\frac{1}{2} u_R(0,t)\Big)^{1-p} 
	= \frac{2^{p-1}}{pt} u_R^{1-p}(0,t)
	\qquad \mbox{for all } r\in (0,r_0),
  \eas
  because $p\ge 1$.
  Upon two integrations using that $\partial_r u_R(0,t)=0$, this first implies that
  \bas
	\partial_r u_R(r,t) \ge -\frac{2^{p-1}}{pt} u_R^{1-p}(0,t) r
	\qquad \mbox{for all } r\in (0,r_0)
  \eas
  and thereafter yields
  \bas
	u_R(r,t) \ge u_R(0,t) - \frac{2^{p-1}}{pt} u_R^{1-p}(0,t)  \frac{r^2}{2}
	\qquad \mbox{for all } r\in [0,r_0].
  \eas
  When evaluated at $r=r_0$, this shows that
  \bas
	\frac{1}{2} u_R(0,t)
	\le \frac{2^{p-2}}{pt} u_R^{1-p}(0,t)  r_0^2
  \eas
  or, equivalently, 
  \be{89.3}
	r_0 \ge \Big(\frac{pt}{2^{p-1}}\Big)^\frac{1}{2} u_R^\frac{p}{2}(0,t).
  \ee
  Since from the definition (\ref{89.2}) of $r_0$ we see that
  \bas
	\int_{B_R} u_R^q(\cdot,t)
	\ge \int_{B_{r_0}} u_R^q(\cdot,t)
	\ge 2^{-q} u_R^q(0,t)  |B_{r_0}|
	= 2^{-q} u_R^q(0,t)  \frac{\omega_n r_0^n}{n},
  \eas
  the inequality (\ref{89.3}) thus entails that
  \bas
	\int_{B_R} u^q(\cdot,t)
	\ge \frac{\omega_n}{2^q n} u_R^q(0,t)  \Big(\frac{pt}{2^{p-1}}\Big)^\frac{n}{2} u_R^\frac{np}{2}(0,t) 
	= \frac{p^\frac{n}{2}\omega_n}{2^{q+\frac{n(p-1)}{2}}n}  t^\frac{n}{2} u_R^\frac{np+2q}{2}(0,t),
  \eas
  which precisely yields (\ref{89.1}), for $u_R(0,t)=\|u_R(\cdot,t)\|_{L^\infty(B_R)}$ again by monotonicity.
\qed
In conjunction with with Lemma~\ref{lem12}, this entails our main result concerning upper estimates for decay 
with respect to the norm in $L^\infty(\R^n)$ of radial and radially nonincreasing solutions emanating from rapidly
decreasing initial data.\abs
\proofc of Theorem~\ref{theo14}.\quad
  We first apply Lemma~\ref{lem12} to find $q>0$, $t_0>0$ and $c_1>0$ such that for any $R>0$,
  the solution of (\ref{0R}) from Lemma \ref{lem202} satisfies
  \be{14.4}
	\|u_R(\cdot,t)\|_{L^q(B_R)}
	\le c_1 t^{-\frac{1}{p}} {\cal L}^{-\frac{np+2q}{npq}} \Big(\frac{1}{t}\Big)
	\qquad \mbox{for all } t\ge t_0.
  \ee
  Thereafter, thanks to the symmetry and monotonicity properties of $u_0$	
  we may invoke Lemma~\ref{lem89} to obtain
  $c_2>0$ fulfilling
  \bas
	\|u_R(\cdot,t)\|_{L^\infty(B_R)} 
	\le c_2 t^{-\frac{n}{np+2q}} \|u_R(\cdot,t)\|_{L^q(B_R)}^\frac{2q}{np+2q}
	\qquad \mbox{for all } t>0.
  \eas
  Combining this with (\ref{14.4}) shows that
  \bas
	\|u_R(\cdot,t)\|_{L^\infty(B_R)}
	&\le& c_2 t^{-\frac{n}{np+2q}}  
	\bigg\{ c_1 t^{-\frac{1}{p}} {\cal L}^{-\frac{np+2q}{npq}} \Big(\frac{1}{t}\Big) \bigg\}^\frac{2q}{np+2q} \\
	&=& c_1^\frac{2q}{np+2q} c_2 t^{-\frac{n}{np+2q} - \frac{2q}{p(np+2q)}} {\cal L}^{-\frac{2}{np}} \Big(\frac{1}{t}\Big) \\
	&=& c_1^\frac{2q}{np+2q} c_2 t^{-\frac{1}{p}} {\cal L}^{-\frac{2}{np}}\Big(\frac{1}{t}\Big)
	\qquad \mbox{for all } t\ge t_0,
  \eas
  which on an application of Fatou's lemma, relying on the approximation properties asserted by Lemma~\ref{lem202}, 
  implies (\ref{14.3}) if we let
  $C:=c_1^\frac{2q}{np+2q} c_2$.
\qed
\subsection{Upper bounds in $L^\infty$: examples}
We next intend to derive Corollary~\ref{cor17} and Corollary~\ref{cor23} by applying Theorem~\ref{theo14}
in the concrete contexts made up by the choices in (\ref{el}).\abs
First concentrating on the former example therein, let us make sure that upon an appropriate 
and essentially trivial extension,
the precise form of the logarithmically fast growth is indeed compatible with both (H)
and the requirements from Section~\ref{sect3.2}.
\begin{lem}\label{lem15}
  Let $\kappa>0, M\ge 2$ and
  \bas
	{\cal L}(s):=\left\{ \begin{array}{ll}
	0, & s=0, \\[1mm]
	\ln^{-\kappa} \frac{M}{s}, \qquad & s\in (0,\frac{M}{2}), \\[1mm]
	\ln^{-\kappa} 2, & s\ge \frac{M}{2}.
	\end{array} \right.
  \eas
  Then ${\cal L}\in C^0([0,\infty)) \cap C^2((0,\frac{M}{2}))$ is positive and nondecreasing on $(0,\infty)$ with
  \be{15.1}
	\frac{d}{ds} \Big(s{\cal L}'(s) \Big) \ge 0
	\qquad \mbox{for all } s\in \Big(0,\frac{M}{2}\Big).
  \ee
  Moreover, given any $\lambda_0>0$ we have
  \be{15.3}
	{\cal L}(s) \le (1+a\lambda) {\cal L}(s^{1+\lambda})
	\qquad \mbox{for all $s>0$ and } \lambda \in (0,\lambda_0),
  \ee
  where
  \be{15.4}
	a:=\left\{ \begin{array}{ll}
	\kappa & \mbox{if } \kappa \le 1, \\[1mm]
	\frac{(1+\lambda_0)^\kappa-1}{\lambda_0} \qquad & \mbox{if } \kappa>1.
	\end{array} \right.
  \ee
\end{lem}
\proof
  To verify (\ref{15.1}), we compute
  \bas
	{\cal L}'(s)=\frac{\kappa}{s} \ln^{-\kappa-1} \frac{M}{s}
	\quad \mbox{and} \quad
	{\cal L}''(s)=-\frac{\kappa}{s^2} \ln^{-\kappa-1} \frac{M}{s}
	+ \frac{\kappa(\kappa+1)}{s^2} \ln^{-\kappa-2} \frac{M}{s},
	\qquad s\in \Big(0,\frac{M}{2}\Big),
  \eas
  whence by positivity of $\kappa$ we indeed obtain that ${\cal L}'(s)>0$ for $s\in (0,\frac{M}{2})$ and that
  \bas
	\frac{s{\cal L}''(s)}{{\cal L}'(s)} 
	\ge \frac{-s\frac{\kappa}{s^2} \ln^{-\kappa-1} \frac{M}{s}}{\frac{\kappa}{s} \ln^{-\kappa-1} \frac{M}{s}}
	= -1
	\qquad \mbox{for all } s\in \Big(0,\frac{M}{2}\Big),
  \eas
  which readily implies (\ref{15.1}).\\
  In proving (\ref{15.3}) we first observe that since $M\ge 2$, for each $s\ge \frac{M}{2}\ge 1$ we trivially have 
  ${\cal L}(s^{1+\lambda}) \ge {\cal L}(s)$ by monotonicity of ${\cal L}$. 
  This implies that we only need to consider the case $s<1$, in which
  $s<\frac{M}{2}$ and also $s^{1+\lambda}<\frac{M}{2}$, so that since 
  $\frac{M}{s^{1+\lambda}}\le \big(\frac{M}{s}\big)^{1+\lambda}$ for all $s>0$ and $\lambda>0$ due to the fact that
  $M\ge 1$, we can estimate
  \be{15.5}
	\frac{{\cal L}(s)}{{\cal L}(s^{1+\lambda})}
	= \frac{\ln^\kappa \frac{M}{s^{1+\lambda}}}{\ln^\kappa \frac{M}{s}}
	\le \frac{\ln^\kappa \big(\frac{M}{s}\big)^{1+\lambda}}{\ln^\kappa \frac{M}{s}}
	= (1+\lambda)^\kappa
	\qquad \mbox{for all } \lambda>0.
  \ee
  Here if $\kappa\le 1$ then by convexity of $0\le \xi\mapsto \xi^\frac{1}{\kappa}$ we have
  $(1+a\lambda)^\frac{1}{\kappa} \ge 1+\frac{a\lambda}{\kappa}=1+\lambda$ according to (\ref{15.4}),
  so that (\ref{15.5}) entails (\ref{15.3}) in this case.\\
  When $\kappa<1$, noting that with $a$ as in (\ref{15.4}), $\psi(\lambda):=1+a\lambda-(1+\lambda)^\kappa$, $\lambda\ge 0$,
  satisfies $\psi(0)=\psi(\lambda_0)=0$ and $\psi''(\lambda)=-\kappa(\kappa-1)(1+\lambda)^{\kappa-2} \le 0$ for all
  $\lambda \ge 0$, we see that $\psi(\lambda)\ge 0$ for all $\lambda\in (0,\lambda_0)$, which combined with
  (\ref{15.5}) completes the proof of (\ref{15.3}).
\qed
A straighforward application of Theorem~\ref{theo14} thus yields the following decay result involving
a precise logarithmic correction to the asymptotics described in Theorem~\ref{theo200} when an appropriate
assumption on fast decay of $u_0$ is formulated as an integrability condition.
\begin{cor}\label{cor16}
  Suppose that $u_0\in C^0(\R^n)$ is positive, radially symmetric and nondecreasing with
  respect to $|x|$ with $u_0(x)\to 0$ as $|x|\to\infty$ and
  \be{16.1}
        \int_{\{u_0<\frac{1}{2}\}} \ln^{-\kappa} \frac{1}{u_0(x)} dx < \infty
  \ee
  for some $\kappa>0$. Then there exist $t_0>1$ and $C>0$ such that the minimal solution of {\rm (\ref{0})} satisfies
  \be{16.2}
        \|u(\cdot,t)\|_{L^\infty(\R^n)}
        \le C t^{-\frac{1}{p}} \ln^\frac{2\kappa}{np} t
        \qquad \mbox{for all } t\ge t_0.
  \ee
\end{cor}
\proof
  Since $u_0$ is bounded, we may choose $M\ge 2$ such that $u_0<\frac{M}{2}$ in $\R^n$, and thereupon let ${\cal L}$ be as
  defined in Lemma~\ref{lem15}. Then using that $M\ge 1$ and that ${\cal L}$ is nondecreasing, we can estimate
  \bas
	\int_{\R^n} {\cal L}(u_0)
	&=& \int_{\{u_0<\frac{1}{2}\}} \ln^{-\kappa} \frac{M}{u_0}
	+ \int_{\{u_0\ge \frac{1}{2}\}} {\cal L}(u_0) \\
	&\le& \int_{\{u_0<\frac{1}{2}\}} \ln^{-\kappa} \frac{1}{u_0} 
	+ {\cal L}\Big(\frac{M}{2}\Big)  \Big| \Big\{u_0\ge\frac{1}{2}\Big\} \Big| \\[2mm]
	&<& \infty
  \eas
  due to (\ref{16.1}) and the fact that $\{u_0\ge \frac{1}{2}\}$ is bounded according to our assumption on asymptotic
  decay of $u_0$.\\
  Consequently, Theorem~\ref{theo14} provides $t_1>0$ and $c_1>0$ such that
  \be{16.3}
	\|u(\cdot,t)\|_{L^\infty(\R^n)}
	\le c_1 t^{-\frac{1}{p}} {\cal L}^{-\frac{2}{np}} \Big(\frac{1}{t}\Big)
	\qquad \mbox{for all } t\ge t_1,
  \ee
  so that if we pick $t_0>\max\{t_1,M\}$, then in particular $\frac{1}{t_0}<\frac{M}{2}$, so that 
  ${\cal L}(\frac{1}{t}) = \ln^{-\kappa} (Mt)$ for all $t\ge t_0$. Since $t_0>M$ furthermore implies that $\ln (Mt) \le 2\ln t$
  for all $t\ge t_0$, (\ref{16.3}) thus yields
  \bas
	\|u(\cdot,t)\|_{L^\infty(\R^n)}
	\le c_1 t^{-\frac{1}{p}} \ln^\frac{2\kappa}{np} (Mt) 
	\le 2^\frac{2\kappa}{np} c_1 t^{-\frac{1}{p}} \ln^\frac{2\kappa}{np} t
	\qquad \mbox{for all } t\ge t_0
  \eas
  and thereby establishes (\ref{16.2}).
\qed
For initial data with the pointwise exponential decay behavior assumed in (\ref{17.1}), 
on a slight shift in the exponent of the respective logarithmic factor, the above integral condition
can be verified, thus yielding temporal decay as claimed.\abs
\proofc of Corollary~\ref{cor17}. \quad
  We pick $c_1\ge c_0$ such that $c_1>1$, and let
  \bas
	\ou_0(x):= c_1 \, e^{-\alpha |x|^\beta},
	\qquad x\in\R^n,
  \eas 
  as well as $r_0:=\big(\frac{2\ln c_1}{\alpha}\big)^\frac{1}{\beta}$. 
  Then for $|x|\ge r_0$ we can estimate
  \bas
	\frac{\ou_0(x)}{e^{-\frac{\alpha}{2}|x|^\beta}}
	= c_1 \, e^{-\frac{\alpha}{2}|x|^\beta}
	\le c_1 \, e^{-\frac{\alpha}{2} \cdot \frac{\ln c_1}{\alpha}} =1,
  \eas
  so writing $\kappa:=\frac{n}{\beta}+\frac{np\delta}{2}$ we have
  \bas
	\ln^{-\kappa} \frac{1}{\ou_0(x)}
	\le \ln^{-\kappa} \Big( e^{\frac{\alpha}{2}|x|^\beta}\Big)
	= c_2 |x|^{-\beta\kappa}
	\qquad \mbox{for all } x\in\R^n \setminus B_{r_0},
  \eas
  with $c_2:=\big(\frac{2}{\alpha}\big)^\kappa$.
  Therefore, 
  \bas
	\int_{\R^n \setminus B_{r_0}} \ln^{-\kappa} \frac{1}{\ou_0(x)} dx
	\le c_2 \int_{\R^n \setminus B_{r_0}} |x|^{-\beta\kappa} dx < \infty
  \eas
  thanks to the fact that $\beta\kappa=n+\frac{np\beta\delta}{2}>n$. As clearly also
  \bas
	\int_{B_{r_0} \cap \{\ou_0<\frac{1}{2}\}} \ln^{-\kappa} \frac{1}{\ou_0(x)} dx
	\le \ln^{-\kappa} 2  \cdot \Big| \Big\{ \ou_0<\frac{1}{2}\Big\} \Big|
  \eas
  is finite, Corollary~\ref{cor16} becomes applicable so as to yield $t_0>1$ and $c_3>0$ such that the minimal solution
  $\ou$ of (\ref{0}) emanating from $\ou_0$ satisfies
  \bea{17.3}
	\|\ou(\cdot,t)\|_{L^\infty(\R^n)}
	\le c_3 t^{-\frac{1}{p}} \ln^\frac{2\kappa}{np} t 
	= c_3 t^{-\frac{1}{p}} \ln^{\frac{2}{p\beta}+\delta} t
	\qquad \mbox{for all } t\ge t_0
  \eea
  according to our definition of $\kappa$.
  Since a comparison argument (\cite{wiegner}) shows that thanks to (\ref{17.1}) we have
  $u_R\le \ou$ in $B_R\times (0,\infty)$ for all $R>0$ and hence
  $u\le \ou$ in $\R^n\times (0,\infty)$, (\ref{17.3}) entails (\ref{17.2}).
\qed
Next, in addressing the second example announced in (\ref{el}) we can proceed quite similarly,
starting with a corresponding counterpart of Lemma~\ref{lem15}.
\begin{lem}\label{lem21}
  Let $\kappa>0, M>e$ and $s_0\in [1,\frac{M}{e})$. Then 
  \bas
	{\cal L}(s):=\left\{ \begin{array}{ll}
	0, & s=0, \\[1mm]
	\ln^{-\kappa} \ln \frac{M}{s}, \qquad & s\in (0,s_0), \\[1mm]
	\ln^{-\kappa} \ln \frac{M}{s_0}, & s\ge s_0,
	\end{array} \right.
  \eas
  defines a function ${\cal L}\in C^0([0,\infty)) \cap C^2((0,s_0))$ which is positive and nondecreasing on $(0,\infty)$ 
  and satisfies
  \be{21.1}
	\frac{d}{ds} \Big(s{\cal L}'(s)\Big) \ge 0
	\qquad \mbox{for all } s\in (0,s_0).
  \ee
  Furthermore, given any $\lambda_0>0$ one can find $a>0$ such that
  \be{21.2}
	{\cal L}(s) \le (1+a\lambda) {\cal L}(s^{1+\lambda})
	\qquad \mbox{for all $s>0$ and } \lambda\in (0,\lambda_0).
  \ee
\end{lem}
\proof
  Let us first observe that ${\cal L}$ indeed is well-defined, positive and nondecreasing on $(0,\infty)$, because
  the assumption $s_0<\frac{M}{e}$ warrants that $\ln \ln \frac{M}{s}>\ln\ln\frac{M}{s_0}>0$ for all $s\in (0,s_0)$.\\
  Now (\ref{21.1}) follows from the fact that for any $s\in (0,s_0)$ we have
  \bas
	{\cal L}'(s)=\frac{\kappa}{s\ln\frac{M}{s}} \ln^{-\kappa-1} \ln\frac{M}{s}
  \eas
  and thus
  \bas
	{\cal L}''(s)
	&=& - \frac{\kappa}{s^2 \ln\frac{M}{s}} \ln^{-\kappa-1}\ln\frac{M}{s}
	+ \frac{\kappa}{s^2 \ln^2 \frac{M}{s}} \ln^{-\kappa-1} \ln\frac{M}{s}
	+ \frac{\kappa(\kappa+1)}{s^2 \ln^2\frac{M}{s}} \ln^{-\kappa-2} \ln\frac{M}{s} \\
	&\ge& - \frac{\kappa}{s^2 \ln\frac{M}{s}} \ln^{-\kappa-1}\ln\frac{M}{s} 
	= -\frac{{\cal L}'(s)}{s},
  \eas
  and in order to verify (\ref{21.2}), given $\lambda_0>0$ we fix $a>0$ large enough such that
  \bas
	\Big(1+\frac{\ln (1+\lambda)}{c_1}\Big)^\kappa \le 1+a\lambda
	\qquad \mbox{for all } \lambda\in (0,\lambda_0),
  \eas
  where $c_1:=\ln\ln \frac{M}{s_0}$.
  Then since $M\ge 1$ and hence $M\le M^{1+\lambda}$ for all $\lambda>0$, in the case when $s<1$, and hence 
  $\max\{s,s^{1+\lambda}\}<s_0$ for all $\lambda>0$, we can estimate
  \bas
	\frac{{\cal L}(s)}{{\cal L}(s^{1+\lambda})}
	&=& \bigg\{ \frac{\ln \ln \frac{M}{s^{1+\lambda}}}{\ln\ln\frac{M}{s}}\bigg\}^\kappa 
	\le \bigg\{ \frac{\ln\ln\big(\frac{M}{s}\big)^{1+\lambda}}{\ln\ln\frac{M}{s}}\bigg\}^\kappa 
	= \bigg\{ 1 + \frac{\ln (1+\lambda)}{\ln\ln\frac{M}{s}}\bigg\}^\kappa \\
	&\le& \bigg\{ 1 + \frac{\ln (1+\lambda)}{c_1} \bigg\}^\kappa 
	\le 1+a\lambda
	\qquad \mbox{for all } \lambda\in (0,\lambda_0).
  \eas
  As by monotonicity of ${\cal L}$ we again have ${\cal L}(s)\le {\cal L}(s^{1+\lambda})$ for all $s\ge 1$ and $\lambda>0$, this proves 
  (\ref{21.2}).
\qed
By Theorem~\ref{theo14}, this again implies a decay estimate, now involving a doubly logarithmic factor,
under a certain integrability condition requiring adequately fast decay of the data.
\begin{cor}\label{cor22}
  Suppose that $u_0\in C^0(\R^n)$ is positive, radially symmetric and nonincreasing
  with respect to $|x|$ and such that $u_0(x)\to 0$ as $|x|\to\infty$ as well as 
  \be{22.1}
        \int_{\{u_0<\eta\}} \ln^{-\kappa} \ln \frac{1}{u_0(x)} dx < \infty
  \ee
  for some $\kappa>0$ and $\eta>0$.
  Then there exist $t_0>e$ and $C>0$ such that for the minimal solution of
  {\rm (\ref{0})} we have
  \be{22.2}
        \|u(\cdot,t)\|_{L^\infty(\R^n)}
        \le C t^{-\frac{1}{p}} \ln^\frac{2\kappa}{np} \ln t
        \qquad \mbox{for all } t\ge t_0.
  \ee
\end{cor}
\proof
  This can be obtained by straightforward adaptation of the argument from Corollary~\ref{cor16}, relying on 
  Lemma~\ref{lem21} rather than Lemma~\ref{lem15}.
\qed
For initial data with doubly exponential decay as in (\ref{23.1}), this can now be seen to imply (\ref{23.2}).\abs
\proofc of Corollary~\ref{cor23}. \quad
  Guided by the procedure from Corollary~\ref{cor17}, we let
  \bas
	\ou_0(x):=c_0 \, \exp \Big\{-\alpha \exp \big(\beta|x|^\gamma\big)\Big\},
	\qquad x\in\R^n,
  \eas
  and note that since
  \bas
	\frac{\exp \Big\{-\alpha \exp \big(\beta|x|^\gamma\big)\Big\}}
		{\exp\Big\{-\exp\big(\frac{\beta}{2}|x|^\gamma\big)\Big\}}
	\to 0
	\qquad \mbox{as } |x|\to\infty,
  \eas
  we can fix $r_0>0$ such that
  \bas
	\ou_0(x) \le \exp \Big\{-\exp \Big(\frac{\beta}{2}|x|^\gamma\Big)\Big\}
	\qquad \mbox{for all } x\in\R^n \setminus B_{r_0}.
  \eas
  Therefore, if we let $\kappa:=\frac{n}{\gamma}+\frac{np\delta}{2}$ and $c_1:=(\frac{2}{\beta})^\kappa$, then
  \bas
	\ln^{-\kappa} \ln \frac{1}{\ou_0(x)}
	\le \ln^{-\kappa} \bigg\{ \exp \Big(\frac{\beta}{2}|x|^\gamma\Big)\bigg\} 
	= c_1 |x|^{-\gamma\kappa}
	\qquad \mbox{for all } x\in \R^n \setminus B_{r_0}
  \eas
  and hence
  \bas
	\int_{\R^n \setminus B_{r_0}} \ln^{-\kappa} \ln \frac{1}{\ou_0(x)} dx
	\le c_1 \int_{\R^n \setminus B_{r_0}} |x|^{-\gamma\kappa} dx < \infty,
  \eas
  because $\gamma\kappa>n$.
  Since this entails that e.g.
  \bas
	\int_{\{\ou_0<e^{-2}\}} \ln^{-\kappa} \ln \frac{1}{\ou_0(x)} dx < \infty,
  \eas
  Corollary~\ref{cor22} provides $t_0>e$ and $c_2>0$ with the property that for the minimal solution $\ou$ of (\ref{0})
  with $\ou|_{t=0}=\ou_0$ we have
  \bas
	\|\ou(\cdot,t)\|_{L^\infty(\R^n)} \le c_2 t^{-\frac{1}{p}} \ln^\frac{2\kappa}{np} \ln (t)
	\qquad \mbox{for all } t\ge t_0.
  \eas
  Since $u\le \ou$ by comparison and (\ref{23.1}), in view of our definition of $\kappa$ this establishes (\ref{23.2}).
\qed
\subsection{Lower estimates: Proof of Theorem~\ref{theo31}}
In order to see that the above decay estimates are essentially optimal, and that hence the interpolation 
inequality from Theorem~\ref{theo11} can as well not be substantially improved any further,
by means of an independent argument based on comparison with separated solutions we finally derive some 
lower bounds for arbitrary positive classical solutions to (\ref{0}) which actually hold in a pointwise sense
for all $(x,t)\in \R^n\times (0,\infty)$.\abs
To prepare this, let us observe that
if $u$ is any positive classical solution of (\ref{0}) in $\R^n\times (0,\infty)$,
then the function $z$ defined on $\R^n\times [0,\infty)$ by letting
\be{transf}
	z(x,\tau):=(t+1)^\frac{1}{p}  u(x,t),
	\qquad \tau=\ln (t+1),
\ee
is a positive classical solution of
\be{0z}
	\left\{ \begin{array}{l}
	z_\tau=z^p \Delta z + \frac{1}{p} z, \qquad x\in \R^n, \ t>0, \\[1mm]
	z(x,0)=u_0(x), \qquad x\in\R^n.
	\end{array} \right.
\ee
To estimate this function from below,
let us first recall the following scaling property of solutions to the associated steady-state problem in
a ball with variable radius (\cite{fast_growth1}).
\begin{lem}\label{lem_w}
  Let $p\ge 1$. For $R>0$, let $w_R \in C^0(\bar B_R) \cap C^2(B_R)$ denote the positive solution of
  \be{w1}
	\left\{ \begin{array}{l}
	-\Delta w_R = \frac{1}{p} w_R^{1-p}, \qquad x\in B_R, \\[1mm]
	w_R=0, \qquad x\in \partial B_R.
	\end{array} \right.
  \ee
  Then for each $R>0$ we have
  \be{w2}
	w_R(x)=R^\frac{2}{p}  w_1\Big(\frac{x}{R}\Big)
	\qquad \mbox{for all } x\in B_R.
  \ee
\end{lem}
Selecting appropriate representatives of this family as spatial profiles of separated solutions of the 
Dirichlet problem for the PDE in (\ref{0z}) in suitable balls, we can indeed achieve the announced
lower estimate for solutions by comparison.\abs
\proofc of Theorem~\ref{theo31}.\quad
  We fix $c_1>0$ such that $pc_1<1$, and given $\tau>0$ we let
  \bas
	R(\tau):=\Lambda^{-1}(c_1\tau).
  \eas
  Then from (\ref{31.1}) we first obtain that $R(\tau)\to\infty$ as $\tau\to\infty$, whereupon a second application
  of (\ref{31.1}) shows that
  \bas
	\frac{2}{\tau} \ln R(\tau)
	= \frac{2c_1}{\Lambda(R(\tau))}  \ln R(\tau) \to 0
	\qquad \mbox{as } \tau\to\infty,
  \eas  
  which since $pc_1<1$ entails the existence of $c_2>0$ such that
  \be{31.33}
	p\Lambda(R(\tau)) + 2\ln R(\tau) \le \tau + c_2
	\qquad \mbox{for all } \tau>0,
  \ee
  because
  \bas
	\frac{p\Lambda (R(\tau)) + 2\ln R(\tau) - \tau}{\tau} 
	= pc_1 + \frac{2}{\tau} \ln R(\tau) - 1 
	\to  pc_1-1
	\qquad \mbox{as } \tau\to\infty.
  \eas
  We now fix $\tau_0>0$ and let $c_3:=1/\|w_1\|_{L^\infty(B_1)}$, as well as
  \be{31.333}
	\delta(\tau_0):=c_3 R^{-\frac{2}{p}}(\tau_0)  e^{-\Lambda(R(\tau_0))},
  \ee
  and
  \be{31.5}
	y(\tau):=\Big\{\delta^{-p}(\tau_0) e^{-\tau} + 1 - e^{-\tau}\Big\}^{-\frac{1}{p}},
	\qquad \tau\ge 0,
  \ee
  observing that $y$ solves
  \be{31.4}
	\left\{ \begin{array}{l}
	y'(\tau)=\frac{1}{p}y(\tau)-\frac{1}{p} y^{p+1}(\tau), \qquad \tau>0, \\[1mm]
	y(0)=\delta(\tau_0).
	\end{array} \right.
  \ee
  Then taking $z$ from (\ref{transf}) and 
  \bas
	\uz(x,\tau):=y(\tau)  w_{R(\tau_0)}(x),
	\qquad x\in \bar B_{R(\tau_0)}, \ \tau\ge 0,
  \eas
  we see that $\uz(x,\tau)=0<v(x,\tau)$ for all $x\in \partial B_{R(\tau_0)}$ and $\tau\ge 0$, and that according to
  the monotonicity of $\Lambda$ and Lemma~\ref{lem_w},
  \bas
	\frac{z(x,0)}{\uz(x,0)}
	&=& \frac{u_0(x)}{\delta(\tau_0)  w_{R(\tau_0)}(x)} 
	\ge \frac{e^{-\Lambda (|x|)}}{\delta(\tau_0)  R^\frac{2}{p}(\tau_0) w_1\big(\frac{x}{R(\tau_0)}\big)} \\
	&>& \frac{e^{-\Lambda(R(\tau_0))}}{\delta(\tau_0) R^\frac{2}{p} \|w_1\|_{L^\infty(B_1)}} 
	= 1
	\qquad \mbox{for all } x\in B_{R(\tau_0)}
  \eas
  thanks to our definition (\ref{31.333}) of $\delta(\tau_0)$. As furthermore by (\ref{w1}) and (\ref{31.4}),
  \bas
	\uz_\tau-\uz^p \Delta \uz - \frac{1}{p} \uz
	&=& y' w_{R(\tau_0)} - y^{p+1} w_{R(\tau_0)}^p \Delta w_{R(\tau_0)} - \frac{1}{p} y w_{R(\tau_0)} \\
	&=& \Big\{ y' - \frac{1}{p} y + \frac{1}{p} y^{p+1} \Big\}  w_{R(\tau_0)} 
	= 0
	\qquad \mbox{in } B_{R(\tau_0)} \times (0,\infty),
  \eas
  a comparison argument (\cite{wiegner}) shows that $z\ge \uz$ in $B_{R(\tau_0)} \times (0,\infty)$.
  When evaluated at $x=0$ and $\tau=\tau_0$, by (\ref{31.5}) and again by
  Lemma~\ref{lem_w} this in particular implies that
  \bas
	z(0,\tau_0) 
	\ge y(\tau_0)  w_{R(\tau_0)}(0) 
	= \Big\{\delta^{-p}(\tau_0) e^{-\tau_0} + 1 - e^{-\tau_0}\Big\}^{-\frac{1}{p}}
	 R^\frac{2}{p}(\tau_0) w_1(0),
  \eas
  where since (\ref{31.333}) and (\ref{31.33}) ensure that
  \bas
	\delta^{-p}(\tau_0) e^{-\tau_0} + 1 - e^{-\tau_0}
	&\le& \delta^{-p}(\tau_0) e^{-\tau_0} + 1 
	= c_3^{-p} R^2(\tau_0) e^{p\Lambda(R(\tau_0))}  e^{-\tau_0} + 1 \\
	&=& c_3^{-p} e^{p\Lambda (R(\tau_0)) + 2\ln R(\tau_0)-\tau_0} +1 
	\le c_3^{-p} e^{c_2} +1,
  \eas
  this shows that with $c_4:=(c_3^{-p} e^{c_2}+1)^{-\frac{1}{p}}  w_1(0)$ we have
  \bas
	z(0,\tau_0)
	\ge c_4 R^\frac{2}{p}(\tau_0)
	= c_4  \Big\{\Lambda^{-1}(c_1\tau_0)\Big\}^\frac{2}{p}
	\qquad \mbox{for all } \tau_0>0.
  \eas
  Transforming back by means of (\ref{transf}), we thus obtain that for each $t>0$, writing $\tau_0:=\ln (t+1)$ 
  we can estimate
  \bea{31.6}
	u(0,t)
	= (t+1)^{-\frac{1}{p}} z(0,\tau_0) 
	\ge (t+1)^{-\frac{1}{p}}  c_4  \Big\{\Lambda^{-1} \big(c_1\ln (t+1)\big)\Big\}^\frac{2}{p}
	\qquad \mbox{for all } t> 0.
  \eea
  Thus, if we pick $t_0\ge 1$ large enough such that $c_1\ln t_0\in \Lambda^{-1}([0,\infty))$, then since 
  $c_1\ln (t+1)\ge c_1\ln t$ and $(t+1)^{-\frac{1}{p}} \ge 2^{-\frac{1}{p}} t^{-\frac{1}{p}}$ for all $t\ge 1$
  we can readily derive (\ref{31.3}) from (\ref{31.6}).
\qed
\subsection{Lower estimates: examples}
The application of the latter to the specific frameworks of initial data satisfying (\ref{32.1})
and (\ref{33.1}), respectively, is now rather straightforward:\abs
\proofc of Corollary~\ref{cor32}. \quad
  We let
  \bas
	\Lambda(s):=\alpha s^\beta - \ln c_0,
	\qquad s\ge 0,
  \eas
  and then see that $\frac{\Lambda(s)}{\ln s}\to +\infty$ as $s\to\infty$, and that $\Lambda$ is strictly increasing on
  $[0,\infty)$ with
  \bas
	\Lambda^{-1}(\sigma) = \Big\{\frac{1}{\alpha} (\sigma+\ln c_0)\Big\}^\frac{1}{\beta}
	\qquad \mbox{for all } \sigma \ge -\ln c_0.
  \eas
  As (\ref{32.1}) warrants that $u_0(x)\ge e^{-\Lambda(|x|)}$ for all $x\in\R^n$, 
  Theorem~\ref{theo31} applies so as to yield $t_1>1$ and $c_1>0$ such that  
  \bas
	\|u(\cdot,t)\|_{L^\infty(\R^n)}
	\ge c_1 t^{-\frac{1}{p}}  \Big\{\Lambda^{-1} (c_1\ln t)\Big\}^\frac{2}{p} 
	= c_1 t^{-\frac{1}{p}}  \Big\{\frac{1}{\alpha} (c_1\ln t + \ln
	c_0)\Big\}^\frac{2}{p\beta}
	\qquad \mbox{for all } t\ge t_1.
  \eas
  Thus, if we pick $t_0\ge t_1$ large enough fulfilling $\ln c_0\ge -\frac{c_1}{2}\ln t_0$, from this we infer that
  \bas
	\|u(\cdot,t)\|_{L^\infty(\R^n)}
	\ge c_1 t^{-\frac{1}{p}}  \Big\{\frac{c_1}{2\alpha} \ln t \Big\}^\frac{2}{p\beta}
	\qquad \mbox{for all } t\ge t_0,
  \eas
  which immediately establishes (\ref{32.2}).
\qed
\proofc of Corollary~\ref{cor33}. \quad
  Proceeding as in the proof of Corollary~\ref{cor32}, we first observe that 
  \bas
	\Lambda(s):=\alpha \, e^{\beta|s|^\gamma} -\ln c_0,
	\qquad s\ge 0,
  \eas
  defines a strictly increasing function on $[0,\infty)$ satisfying
  $\frac{\Lambda(s)}{\ln s}\to+\infty$ as $s\to\infty$ as well as
  \bas
	\Lambda^{-1}(\sigma) 
	= \bigg\{ \frac{1}{\beta}\ln \Big[ \frac{1}{\alpha} (\sigma+\ln c_0)\Big] \bigg\}^\frac{1}{\gamma}
	\qquad \mbox{for all } \sigma\ge \alpha-\ln c_0.
  \eas
  According to (\ref{33.1}), Theorem~\ref{theo31} therefore provides $t_1>e$ and $c_1>0$ such that
  \bas
	\|u(\cdot,t)\|_{L^\infty(\R^n)} 
	\ge c_1 t^{-\frac{1}{p}}  \Big\{\Lambda^{-1}(c_1 \ln t)\Big\}^\frac{2}{p} 
	= c_1 t^{-\frac{1}{p}} 
	\bigg\{\frac{1}{\beta} \ln \Big[\frac{1}{\alpha} (c_1\ln t + \ln c_0)\Big] \bigg\}^\frac{2}{p\gamma}
	\qquad \mbox{for all }t\ge t_1.
  \eas
  Hence, picking $t_0\ge t_1$ in such a way that $\ln c_0\ge -\frac{c_1}{2}\ln t_0$ and 
  $\ln\frac{c_1}{2\alpha}\ge -\frac{1}{2}\ln\ln t_0$, we conclude that
  \bas
	\|u(\cdot,t)\|_{L^\infty(\R^n)}
	&\ge& c_1 t^{-\frac{1}{p}}  
	\bigg\{ \frac{1}{\beta} \ln \Big[ \frac{c_1}{2\alpha} \ln t \Big] \bigg\}^\frac{2}{p\gamma} 
	= c_1 t^{-\frac{1}{p}}  
	\bigg\{ \frac{1}{\beta} \Big[ \ln \frac{c_1}{2\alpha} + \ln\ln t\Big] \bigg\}^\frac{2}{p\gamma} \\
	&\ge& c_1 t^{-\frac{1}{p}}  \Big\{ \frac{1}{2\beta} \ln\ln t\Big\}^\frac{2}{p\gamma}
	\qquad \mbox{for all } t\ge t_0,
  \eas
  which clearly entails (\ref{33.2}).
\qed

\noindent
{\bf Acknowledgements.} The first author was supported in part by the Slovak
Research and
Development Agency under the contract No. APVV-14-0378 and by the VEGA grant
1/0319/15.
\end{document}